\documentclass[a4paper,12pt,leqno, twoside]{article}
\usepackage{amssymb}
\usepackage{amsthm}
\usepackage{amstext}
\usepackage{amsbsy}
\usepackage{amscd}
\usepackage{amsopn}
\usepackage{amsmath}
\usepackage{color}
\usepackage{mathtools}
\usepackage{mathdots}
\usepackage[mathscr]{eucal}
\usepackage{graphicx}
\usepackage{makeidx}
\makeindex
\usepackage{epic,eepic}
\usepackage{epsfig}
\usepackage[all]{xy}
\usepackage[french]{babel}
\usepackage[latin1]{inputenc}
\usepackage[nottoc]{tocbibind}

    \def\AM{{\mathbb{A}}}

    \def\FM{{\mathbb{F}}}
    
    \def\HM{{\mathbb{H}}}

    \def\NM{{\mathbb{N}}}
    
    \def\PM{{\mathbb{P}}}
    \def\QM{{\mathbb{Q}}}
    
\def\SG{{\mathfrak S}}

\def\XG{{\mathfrak X}}    
    
    \def\ZM{{\mathbb{Z}}}

    \def\BC{{\mathcal{B}}}
    \def\CC{{\mathcal{C}}}

\def\Hb{{\mathbf H}}    \def\HC{{\mathcal{H}}}
    \def\IC{{\mathcal{I}}}

    \def\LC{{\mathcal{L}}}
    \def\MC{{\mathcal{M}}}
    
    \def\OC{{\mathcal{O}}}

  \def\rb{{\mathbf r}}  
    
\def\Tb{{\mathbf T}}    \def\TC{{\mathcal{T}}}
  \def\ub{{\mathbf u}}

\def\a{\alpha}
\def\b{\beta}

\def\G{\Gamma}
\def\d{\delta}
\def\D{\Delta}
\def\e{\varepsilon}
\def\ph{\varphi}
\def\l{\lambda}
\def\L{\Lambda}

\def\o{\omega}
\def\O{\Omega}

\def\s{\sigma}
\def\Sig{\Sigma}
\def\th{\theta}

\def\z{\zeta}

\def\psib{{\boldsymbol{\psi}}}

\def\etab{{\boldsymbol{\eta}}}          

\DeclareMathOperator{\End}{{\mathrm{End}}}

\DeclareMathOperator{\Fr}{{\mathrm{Fr}}}
\DeclareMathOperator{\Gal}{{\mathrm{Gal}}}
\DeclareMathOperator{\Hom}{{\mathrm{Hom}}}

\DeclareMathOperator{\ind}{{\mathrm{ind}}}

\DeclareMathOperator{\Ker}{{\mathrm{Ker}}}

\DeclareMathOperator{\Proj}{{\mathrm{Proj}}}

\DeclareMathOperator{\res}{{\mathrm{res}}}

\DeclareMathOperator{\Sp}{{\mathrm{Sp}}}
\DeclareMathOperator{\Spec}{{\mathrm{Spec}}}

\DeclareMathOperator{\Stab}{{\mathrm{Stab}}}

\DeclareMathOperator{\val}{{\mathrm{val}}}

\def\oQl{\o\QM_{\ell}}

\def\oFq{\o{\FM}_q}

\def\para{sous-groupe parabolique }

\def\Omenr{\wh{\Omega}^{d-1}_{\breve{\OC}}}
\def\Ome{\Omega_{\breve{K}}^{d-1}}

\def\Art {\mathrm {Art}}
\def\Nr {\mathrm {Nr}}

\def\Gal{{\rm Gal}}
\def\ord{{\rm ord}}
\def\Frob{{\rm Frob}}
\def\Nilp{{\rm Nilp}}
\def\Proj{{\rm Proj}}
\def\Sym{{\rm Sym}}
\def\Isom{{\rm Isom}}
\def\Pic{{\rm Pic}}
\def\lie{{\rm Lie}}
\def\can{{\rm can}}
\def\spe{{\rm sp}}
\def\rig{{\rm rig}}

\def\ddonc{\,\,\Longrightarrow\,\,}

\def\simto{\buildrel\hbox{$\sim$}\over\longrightarrow}

\def\leq{\leqslant}
\def\geq{\geqslant}
\def\into{\hookrightarrow}
\def\onto{\twoheadrightarrow}
\def\id{\mathop{\mathrm{Id}}\nolimits}

\def\ba{\backslash}
\def\wt{\widetilde}
\def\wh{\widehat}
\def\o#1{\overline{#1}}

\def\To#1{\buildrel\hbox{\tiny{$#1$}}\over\longrightarrow}

\def\to{\rightarrow}

\def\Hom{\mathop{\hbox{\rm Hom}}\nolimits}

\def\dim{\mathop{\mbox{\rm dim}}\nolimits}

\def \limi#1{\lim\limits_{\displaystyle\longrightarrow\atop {#1}}}

\def\ini{\setcounter{equation}{\value{subsubsection}}\addtocounter{subsubsection}{1}}

\makeatletter
\renewcommand{\subsubsection}{\@startsection{subsubsection}{3}{\parindent}{-\baselineskip}{-0.01\baselineskip}{\bf}}
\renewcommand*{\@seccntformat}[1]{%
  \csname the#1\endcsname\
} \makeatother

\def\ali{\subsubsection{}\setcounter{equation}{0}}

\newtheoremstyle{th}
  {\baselineskip}{.5\baselineskip}{\itshape}
  {\parindent}{\bf}
  {--}{.5em}{}

\newtheoremstyle{def}
  {\baselineskip}{\baselineskip}{}
  {\parindent}{\bf}
  {--}{.5em}{}

\newtheoremstyle{th*}
  {.5\baselineskip}{.5\baselineskip}{\itshape}
  {\parindent}{\bf}
  {--}{.5em}{}

\newtheoremstyle{remark*}
  {.5\baselineskip}{.5\baselineskip}{}
  {\parindent}{\bf}
  {--}{.5em}{}

\newtheoremstyle{remark}
  {.5\baselineskip}{.5\baselineskip}{}
  {\parindent}{\bf}
  {--}{.5em}{}

\swapnumbers \theoremstyle{th}
\newtheorem{theo}[subsubsection]{\sc Th{\'e}or{\`e}me.\bf}
\newtheorem{lemme}[subsubsection]{\sc Lemme.\bf}
\newtheorem{prop}[subsubsection]{\sc Proposition.\bf}
\newtheorem{coro}[subsubsection]{\sc Corollaire.\bf}

\theoremstyle{def}
\newtheorem{fact}[subsubsection]{\sc Fait\bf}

\newtheorem{DEf}[subsubsection]{\sc D{\'e}finition.\bf}

\theoremstyle{remark}
\newtheorem{exem}[subsubsection]{\sc Exemple.\bf}   

\theoremstyle{th*}
\newtheorem*{thm}{\sc Th{\'e}or{\`e}me.}
\newtheorem*{lem}{\sc Lemme.}

\theoremstyle{remark*}
\newtheorem*{rem}{\sc Remarque.}
\newtheorem*{fait}{\sc Fait.}

\newcommand{\findem}{\hfill$\Box$\par\medskip}
\newcommand{\dem}{\indent {\it Preuve :} \rm }

\newenvironment{preuve}{\dem}{\findem}

\title{L'espace sym\'etrique de Drinfeld et correspondance\\
 de Langlands locale I}
\author{Haoran Wang}
\date{}

\theoremstyle{plain}

\begin{document}
\maketitle
\renewcommand{\proofname}{\indent Preuve}

\def\dd{D_d^\times}
\def\mdro{\MC_{Dr,0}}
\def\mdrn{\MC_{Dr,n}}
\def\mdr{\MC_{Dr}}
\def\mlto{\MC_{LT,0}}
\def\mltn{\MC_{LT,n}}
\def\mlt{\MC_{LT}}
\def\mltK{\MC_{LT,K}}
\def\LJ{{\rm LJ}}
\def\JL{{\rm JL}}
\def\SL{{\rm SL}}
\def\GL{{\rm GL}}
\def\Spf{{\rm Spf}}
\def\DL{{\rm DL}}
\def\PGL{{\rm PGL}}
\def\Coef{{\rm Coef}}

\def\Ql{\QM_{\ell}}
\def\Zl{\ZM_{\ell}}
\def\Fl{\FM_{\ell}}

\abstract{Soit $K$ une extension finie de $\QM_p.$ On \'etudie la g\'eom\'etrie et la cohomologie du rev\^etement mod\'er\'e de l'espace sym\'etrique de Drinfeld sur $K.$ On prouve, de mani\`ere purement locale, que la cohomologie de degr\'e m\'edian r\'ealise la correspondance de Langlands locale et la correspondance de Jacquet-Langlands locale pour les repr\'esentations supercuspidales de niveau z\'ero.}

\begin{center}
{\bf {Abstract}}
\end{center}

{Let $K$ be a finite extension of $\QM_p.$ We study the geometry and cohomology of the tamely ramified cover of Drinfeld's symmetric space over $K.$ We prove, in a purely local way, that the cohomology in middle degree realises the local Langlands correspondence and local Jacquet-Langlands correspondence for the depth zero supercuspidal representations.}

\tableofcontents

\section{Introduction}

Soit $K$ un corps local de caract\'eristique r\'esiduelle $p,$ d'anneau des entiers $\OC$ et de corps r\'esiduel $\FM_q.$ Fixons une cl\^oture alg\'ebrique $K^{ca}$ de $K$ et un entier $d\geq 2.$ On dispose de trois groupes $G:=\GL_d(K),$ $D^\times$ les \'el\'ements inversibles dans l'alg\`ebre \`a division $D$ qui est centrale sur $K$ et est d'invariant $1/d,$ et $W_K$ le groupe de Weil de $K.$

Dans \cite{drinfeld-ell}, Drinfeld a introduit l'``espace sym\'etrique $p$-adique'' de dimension $d-1,$ $\O^{d-1},$ d\'efini comme le compl\'ementaire de tous les hyperplans $K$-rationels dans l'espace projectif $\PM^{d-1}_K.$ Cet espace est un espace rigide-analytique muni d'une action continue de $G.$ Peu apr\`es, il a d\'ecouvert, dans \cite{drinfeld-covering}, un syst\`eme projectif de rev\^etements \'etales $\{\Sig_n\}$ (la tour de Drinfeld) de $\O^{d-1}$ munis des actions de $G\times D^\times\times W_K.$ Drinfeld et Carayol \cite{Carayol-LT} ont conjectur\'e que la limite inductive de la cohomologie de $\Sig_n$ se d\'ecompose en une somme directe sur les s\'eries discr\`etes $\pi$ de $G,$ de $\pi\otimes \LJ(\pi)\otimes \s(\pi),$ o\`u $\LJ(\pi)$ d\'esigne la repr\'esentation de $D^\times$ associ\'ee \`a $\pi$ par la correspondance de Jacquet-Langlands, et $\s(\pi)$ est l'unique quotient irr\'eductible de la repr\'esentation de Weil-Deligne $L(\pi)$ associ\'ee \`a $\pi$ par la correspondance de Langlands normalis\'ee \`a la Hecke. Lorsque $\pi$ est supercuspidale, $\s(\pi)$ co\"incide avec $L(\pi).$

La partie {\em supercuspidale} de cette conjecture ainsi que l'\'enonc\'e analogue concernant la tour de Lubin-Tate $\{\MC_{LT,n}\}_{n\in \NM}$ (la conjecture de Deligne-Carayol) ont \'et\'e explor\'es dans les quatres articles \cite{Harris-Carayol} \cite{Boyer-these} \cite{Harris-Taylor} et \cite{Hausberger-these}, en caract\'eristique nulle ou en \'egales caract\'eristique, par voie \emph{globale} en utilisant les vari\'et\'es de Shimura ou de Drinfeld.

Comme la th\'eorie de Lubin-Tate classique, il est naturel de chercher une preuve de voie {\em locale} et plus explicite de ces r\'esultats. L'objet de cet article est d'\'etudier par voie {\em locale} la partie {\em supercuspidale} de la cohomologie de $\Sig_1$ (par abus de notation, on le note $\Sig$ dans les sections qui suivent) lorsque $K$ est de caract\'eristique nulle (cette \'etude est aussi valable pour $K$ d'\'egale caract\'eristique). La partie {\em non-supercuspidale} sera l'objet de \cite{Wang-Sigma2}. Pour \'enoncer nos r\'esultats, on utilise une variante de $\Sig_1$ not\'ee $\MC_{Dr,1}/\varpi^\ZM$ qui est une r\'eunion disjointe de $d$-copies de $\Sig_1.$ Le th\'eor\`eme principal est le suivant:

\medskip
\noindent{\bf Th\'eor\`eme A.}
(Th\'eor\`eme \ref{T:1})
{\it
Soit $\rho$ une $\oQl$-repr\'esentation irr\'eductible de niveau z\'ero de $D^\times$ de caract\`ere central trivial sur $\varpi^\ZM$ telle que $\JL(\rho)$ soit une repr\'esentation supercuspidale de $G$ par la correspondance de Jacquet-Langlands. Alors, en tant que representations de $G\times W_K,$ on a
      $$
\Hom_{D^\times}\big(\rho,H^{i}_c(\MC_{Dr,1}/\varpi^\ZM,\oQl)\big)\cong\left\{
                                       \begin{array}{ll}
                                         \JL(\rho)\otimes L(\JL(\rho)), & \hbox{si $i=d-1$;} \\
                                         0, & \hbox{sinon.}
                                       \end{array}
                                     \right.
$$
}

\medskip

Notre preuve repose sur l'\'etude de la g\'eom\'etrie de $\Sig_1.$ Rappelons qu'il existe une application $\tau:\O^{d-1}\To{}|\BC\TC|,$ o\`u $|\BC\TC|$ d\'esigne la r\'ealisation g\'eom\'etrique de l'immeuble de Bruhat-Tits semi-simple $\BC\TC$ associ\'e \`a $G.$ En prenant la composition avec la transition $p:\Sig_1\to\O^{d-1}_K,$ on obtient un morphisme $\nu:\Sig_1\To{}|\BC\TC|.$ Comme il n'existe pas de notion de ``base de Drinfeld'' \`a l'instant, on construit sur chaque sommet $s\in\BC\TC$ un mod\`ele entier {\em lisse} de $\nu^{-1}(|s|)$ dans \ref{Sec:4}, et on obtient en particulier le r\'esultat suivant:


\medskip

\noindent{\bf Th\'eor\`eme B.} (Th\'eor\`eme \ref{C1C1})
{\it Si $s$ est un sommet de $\BC\TC,$ on a un isomorphisme:
$$
H^{q}_c(\nu^{-1}(|s|^*),\L)\simto H^q_c(\DL^{d-1}_{\oFq},\L),
$$
o\`u $|s|^*=\cup_{s\in\s}|\s|$ et $\DL^{d-1}_{\oFq}$ d\'esigne le rev\^etement de Deligne-Lusztig Coxeter associ\'e \`a $\GL_d(\FM_q),$}

\bigskip

Ce travail est inspir\'e par celui de Genestier \cite{genestier} sur l'espace de Drinfeld mod\'er\'e lorsque $K$ est de caract\'eristique positive, et par celui de Yoshida \cite{Yoshida} sur le niveau mod\'er\'e de la tour de Lubin-Tate. Notre r\'esultat est n\'eanmoins plus pr\'ecis que celui de Yoshida qui ne d\'ecrit que l'action de $G$ et du groupe d'inertie. Notre d\'emonstration est inspir\'ee par Teitelbaum \cite{teitelbaum} qui \'etudie la g\'eom\'etrie de $\Sig_1$ lorsque $d=2.$ Dans \cite{Strauch}, Strauch a \'etudi\'e la correspondance de Jacquet-Langlands dans la tour de Lubin-Tate.

Nous d\'ecrivons bri\`evement les contenus de diff\'erents paragraphes. Dans les paragraphes 2.1 et 2.2, on rappelle l'espace de Drinfeld $\O^{d-1}$ et son rev\^etement mod\'er\'e $\Sig_1.$ On d\'emontre au paragraphe 2.2, en employant un lemme de Zheng \cite{zheng}, que l'\'etude de la cohomologie sans support $H^q(\nu^{-1}(|s|^*),\L)$ se ram\`ene \`a celle de $H^q(\nu^{-1}(|s|),\L).$ Au paragraphe \ref{S2Pa3}, on construit un mod\`ele entier {\em lisse} de $\nu^{-1}(|s|),$ d\'efini sur une extension mod\'er\'ement ramifi\'ee du compl\'et\'e de l'extension non ramifi\'ee maximale de $K,$ dont la fibre sp\'eciale $\o\Sig^0_s$ calcule la cohomologie sans support de $\nu^{-1}(|s|).$ Au cours de la preuve, on utilise la classification, d\^ue \`a Raynaud \cite{Raynaud}, des sch\'emas en groupes de type $(p,\ldots,p).$ Aux paragraphes 2.4 et 2.5, on d\'emontre que $\o\Sig^0_s\cong \DL^{d-1}_{\oFq}$ (et en cons\'equence le Th\'eor\`eme B.), en calculant leurs classes de $\mu_{q^d-1}$-torseur dans $H^1(\O^{d-1}_{\oFq},\mu_{q^d-1}),$ o\`u $\O^{d-1}_{\oFq}$ d\'esigne l'espace de Drinfeld sur le corps fini $\FM_q.$ Le th\'eor\`eme A. est alors obtenu dans la section 3. Pour r\'ealiser la correspondance de Langlands locale dans la cohomologie, on a besoin d'un r\'esultat sur les composantes connexes de $\Sig_1.$ Ce r\'esultat \'etant connu dans le cas d'\'egales caract\'eristiques \cite{genestier-livre}, peut \^etre obtenu par le travail de Chen \cite{chen} sur le c\^ot\'e de Lubin-Tate, en vertu d'isomorphisme de Faltings-Fargues, et il nous permet de descendre la fibre sp\'eciale $\o\Sig^0_s,$ et donc de d\'ecrire l'action de Frobenius.

\medskip
\textbf{Remerciements:} Je remercie profond\'ement mon directeur de th\`ese Jean-Fran\c cois Dat de me proposer ce sujet et ses constants encouragements. Je tiens \`a exprimer ma gratitude \`a Miaofen Chen, Bas Edixhoven, Jared Weinstein et Weizhe Zheng pour leurs enrichissantes discussions. Je remercie Jared Weinstein pour son int\'er\^et sur ce travail et pour m'avoir invit\'e \`a Boston. Enfin, je remercie le referee anonyme pour ses commentaires.   

\section{Sur le rev\^etement mod\'er\'e de l'espace sym\'etrique de Drinfeld}

Dans cette section, on rappelle tout d'abord l'espace sym\'etrique $p$-adique de Drinfeld introduit dans \cite{drinfeld-ell}. C'est un espace rigide-analytique dont un mod\`ele entier param\`etre certains groupes $p$-divisibles \cite{drinfeld-covering}. On rappelle le niveau mod\'er\'e de la tour de Drinfeld. Le lecteur pourra consulter \cite{boutot-carayol} pour une d\'emontration d\'etaill\'ee de l'interpretation modulaire lorsque $d=2.$ Ensuite, on donne des r\'esultats sur la g\'eom\'etrie du niveau mod\'er\'e.

\subsection{Rappels sur l'espace sym\'etrique de Drinfeld\underline{\underline{}}}\label{ee}

\ali\label{Sec:5} Soit $d\geq2$ un entier. Fixons un corps $p$-adique $K,$ une extension finie de
$\QM_p,$ et notons $\OC$ son anneau des entiers et $\varpi$ une uniformisante de
$\OC$. Le corps r\'esiduel $\OC/\varpi\simeq\FM_q$ est une extension
finie de $\FM_p$, son degr\'e $[\FM_q:\FM_p]$ sera not\'e $f$. On
fixe $K^{ca}$ une cl\^oture alg\'ebrique de $K$ et $\widehat{K^{ca}}$
son compl\'et\'e par l'unique norme qui \'etend celle de $K.$ Soient $D$ l'alg\`ebre \`a division centrale sur
$K$ d'invariant $1/d$ et $\OC_D$ l'anneau des entiers de $D$. Soient
$K_d$ une extension non-ramifi\'ee de degr\'e $d$ de $K$ contenue
dans $D$, $\OC_d$ l'anneau des entiers de $K_d$. Il existe un
\'el\'ement $\Pi_D\in \OC_D$ tel que $\OC_D$ soit engendr\'e sur
$\OC_d$ par $\Pi_D$ v\'erifiant les relations $\Pi^d_D=\varpi$ et $\Pi_D
a=\sigma(a)\Pi_D$, $\forall a\in \OC_d$, o\`u $\sigma\in\Gal(K_d/K)$
le rel\`evement de Frobenius. On note $\breve{K}$ le compl\'et\'e de
l'extension non-ramifi\'ee maximale de $K$ dans la cl\^oture
alg\'ebrique $K^{ca},$ $\breve{\OC}$ son anneau des entiers, on a donc
$\breve{K}=\breve{\OC}[1/p]$.

Notons $\BC\TC$ l'immeuble de Bruhat-Tits semi-simple de $G:=\GL_d(K),$ c'est un complexe simplicial dont l'ensemble de sommets s'identifie \`a l'ensemble des classes d'homoth\'etie de $\OC$-r\'eseaux dans l'espace vectoriel $K^d.$ Un ensemble de sommets $\s=\{s_0,\ldots,s_k\}$ forme un $k$-simplexe s'il existe des repr\'esentants $\L_i$ de $s_i~\forall 0\leq i\leq k$ tels que $\varpi \Lambda_k\subsetneq\Lambda_0\subsetneq\cdots\subsetneq\Lambda_{k}.$ On notera $\BC\TC_{k}$ l'ensemble des $k$-simplexes. On d\'esignera $|\BC\TC|$ la r\'ealisation g\'eom\'etrique de $\BC\TC.$ Pour $\s=\{\varpi
\Lambda_k\subsetneq\Lambda_0\subsetneq\cdots\subsetneq\Lambda_{k}\}$
un $k$-simplexe de $\BC\TC,$ notons $|\s|\subset|\BC\TC|$ sa facette associ\'ee: on a pour toute paire de simplexes $\s'\subset\s,$ $|\s'|\subset\o{|\s|},$ o\`u $\o{|\s|}$ d\'esigne l'adh\'erence de $|\s|.$ On notera aussi $|\s|^*\subset|\BC\TC|$ la r\'eunion de toutes les facettes $|\s'|$ avec $\s'$ contenant $\s$, i.e. $|\s|^*=\bigcup_{\s\subset\s'}|\s'|.$ \'Evidemment, si $\s=\{s_0,\ldots,s_k\}$ avec $s_i\in\BC\TC_{0}$ des sommets, on a $|\s|^*=\bigcap_i|s_i|^*.$

Notons $\Omega_K^{d-1}$ l'espace sym\'etrique de Drinfeld de dimension $d-1,$ d\'efini dans \cite{drinfeld-ell} comme un sous $K$-espace rigide-analytique de l'espace projectif $\PM^{d-1}_K.$ Il s'identifie au
compl\'ementaire de l'ensemble des hyperplans $K$-rationnels dans $\PM^{d-1}_K.$ 
On sait que les points de $|\BC\TC|$ s'identifient aux classes d'homoth\'etie de {\em normes} sur le $K$-espace vectoriel $K^d$ ({\em cf.} \cite{GI}). Ceci nous fournit une {\em application de r\'eduction} $\tau:\Omega_K^{d-1}\to|\BC\TC|$ (voir \cite{DH}). On d\'esignera par le m\^eme symbole $\Omega_K^{d-1}$ le $K$-espace analytique {\em \`a la} Berkovich correspondant, et il est naturellement muni d'une action \emph{continue} de $G$ triviale sur le centre ({\em cf.} \cite{berkovich-drin}). On obtient alors un espace analytique $\Omega_K^{d-1,ca}:=\Omega_K^{d-1}\widehat{\otimes}_K\widehat{K^{ca}}$ muni d'une action continue du groupe de Weil $W_K$ de $K$ en \'etendant les scalaires \`a $\widehat{K^{ca}}.$

\ali Dans un travail non publi\'e, Deligne introduit un mod\`ele
semi-stable $\wh{\Omega}^{d-1}_\OC$ de $\Omega_K^{d-1}$ sur $\Spf \OC,$ en recollant les mod\`eles locaux (au-dessus des simplexes de $\BC\TC$) que nous rappelons ci-dessous ({\em cf.} \cite{Rapoport}).

On note alors $\wh{\Omega}_{\OC,\s}^{d-1}$ le sch\'ema formel classifiant les classes d'isomorphie de
diagrammes commutatifs
$$
\xymatrix{
  \Lambda_{-1}=\varpi\Lambda_k \ar[d]_{\alpha_k/\varpi} \ar@{^(->}[r] & \Lambda_0 \ar[d]_{\alpha_0} \ar@{^(->}[r] & \dots
{}
  \ar@{^(->}[r] & \Lambda_k \ar[d]^{\alpha_k} \\
  L_k\ar[r]^{\Pi} \ar@(rd,ld)[rrr]_{\times\varpi}  & L_0 \ar[r]^{\Pi} & \dots \ar[r]^{\Pi} & L_k   }
$$
sur un $\OC$-sch\'ema $S=\Spec(R)$, o\`u $R\in \Nilp_\OC$, la cat\'egorie
des $\OC$-alg\`ebres sur lesquelles $\varpi$ est nilpotent. Les $L_i$
sont des fibr\'es en droites sur $S$, les applications $\alpha_i$
sont des homomorphismes de $\OC$-modules, les morphismes $\Pi$ sont
des homomorphismes de $\OC_S$-modules, v\'erifiant la condition:
pour tout $x\in\Spec(R)$, on a
$$\Ker\{\alpha_i(x):\Lambda_i/ \varpi\Lambda_i\to L_i\otimes_R
k(x)\}\subset\Lambda_{i-1}/\varpi\Lambda_i,$$
o\`u $k(x)$ est le corps r\'esiduel de $x.$ L'objet universel sur $\wh{\Omega}_{\OC,\s}^{d-1}$ sera not\'e:
$$
\xymatrix{
  \varpi\Lambda_k \ar[d]_{\alpha_k/\varpi} \ar@{^(->}[r] & \Lambda_0 \ar[d]_{\alpha_0} \ar@{^(->}[r] & \dots
{}
  \ar@{^(->}[r] & \Lambda_k \ar[d]^{\alpha_k} \\
  \LC_k \ar[r]^{\Pi} & \LC_0 \ar[r]^{\Pi} & \dots \ar[r]^{\Pi} & \LC_k   }
$$

\begin{fait}({\em cf.} \cite[(III.1)]{genestier-livre})
On note $\wt{\PM}_\s$ le $\OC$-sch\'ema obtenu \`a partir de
$\PM(\Lambda_k)=\Proj(\Sym(\Lambda_k))$ en l'\'eclatant successivement
le long du sous-sch\'ema ferm\'e $\PM(\Lambda_k/\Lambda_{k-1})$ de
sa fibre sp\'eciale
$\wt{\PM}_{\s}\otimes\oFq=\PM(\Lambda_k/\varpi\Lambda_k)$, puis en
\'eclatant le transform\'e strict de $\PM(\Lambda_k/\Lambda_{k-2})$,
puis le transform\'e strict de $\PM(\Lambda_k/\Lambda_{k-3})$ et
ainsi de suite. Comme expliqu\'e dans \cite[(III.1)]{genestier-livre}, $\wt{\PM}_\s$ repr\'esente un foncteur $H_\s$ sur $\Nilp_\OC$ qui associe \`a $R\in\Nilp_\OC$ l'ensemble des classes d'isomorphie de diagrammes commutatifs

$$
\xymatrix{
  \Lambda_{-1}=\varpi\Lambda_k \ar[d]_{\alpha_k/\varpi} \ar@{^(->}[r] & \Lambda_0 \ar[d]_{\alpha_0} \ar@{^(->}[r] & \dots
{}
  \ar@{^(->}[r] & \Lambda_k \ar[d]^{\alpha_k} \\
  L_k\ar[r]^{\Pi} \ar@(rd,ld)[rrr]_{\times\varpi}  & L_0 \ar[r]^{\Pi} & \dots \ar[r]^{\Pi} & L_k   }
$$
o\`u les $L_i$ sont des $R$-modules inversibles, les applications $\alpha_i$
sont des homomorphismes de $\OC$-modules, les morphismes $\Pi$ sont
des homomorphismes de $R$-modules, tels que pour tout $i$ le morphisme $B$-lin\'eaire 
$$
\a_i\otimes \id_R:\L_i\otimes_{\OC}R\To{}L_i
$$
soit surjectif. On a alors une immersion de $\wh{\Omega}^{d-1}_{\OC,\s}$ dans $\wt{\PM}_\s.$ On d\'esigne $\Omega^0_{\OC,\s}$ le compl\'ementaire dans la fibre sp\'eciale $\wt{\PM}_{\s}\otimes\oFq$ de $\wt{\PM}_{\s}$ des ferm\'es
\[
\{\alpha_j(m)=0\}\  \ (j\in\ZM/k\ZM,
m\in\Lambda_j/\varpi\Lambda_j-\Lambda_{j-1}/\varpi\Lambda_j).
\]
Alors $\wh{\Omega}^{d-1}_{\OC,\s}$ s'identifie au compl\'et\'e de
$\wt{\PM}_{\s}$ le long de $\Omega^0_{\OC,\s}.$  
Lorsque $\s'$ est un sous-simplexe de $\s$, on a une immersion naturelle ouverte $\wh{\Omega}^{d-1}_{\OC,\s'}\into \wh{\Omega}^{d-1}_{\OC,\s}.$ Le sch\'ema formel $\wh{\Omega}^{d-1}_\OC$ est par d\'efinition la limite inductive $\limi{\s\in \BC\TC}\wh{\Omega}^{d-1}_{\OC,\s}.$ De plus, la fibre g\'en\'erique de $\wh{\Omega}^{d-1}_{\OC}$ (au sens de Raynaud-Berkovich) s'identifie \`a $\O^{d-1}_K.$ Si $g\in G$ envoyant un simplexe $\s$ sur $g\s,$ on d\'efinit un isomorphisme $g:\wh{\Omega}^{d-1}_{\OC,\s}\to \wh{\Omega}^{d-1}_{\OC,g\s}$ en associant \`a une donn\'ee $(\a_j:\L_j\to L_j,\Pi)$ la donn\'ee $(g\L_j\To{g^{-1}}\L_j\To{\a_j}L_j,\Pi).$ Ceci d\'efinit une action de $G$ sur le syst\`eme inductif $(\wh{\Omega}^{d-1}_{\OC,\s})_\s$, et donc sur $\O^{d-1}_K.$ Sous cette action, l'application de r\'eduction $\tau:\O^{d-1}_K\to|\BC\TC|$ est $G$-\'equivariante.
\end{fait}

\begin{DEf}\label{Def::1}
Soit $\s=\{\varpi
\Lambda_k\subsetneq\Lambda_0\subsetneq\cdots\subsetneq\Lambda_{k}\}$
un $k$-simplexe, le {\em type} de $\s$ est la suite des entiers
$(e_0,\ldots,e_k)$ telle que
$e_0=\dim_{\FM_q}\Lambda_0/\varpi\Lambda_k$ et
$e_i=\dim_{\FM_q}\Lambda_i/\Lambda_{i-1}$ pour $1\leq i\leq k$.
\'Evidemment, on a $\sum e_i=d$.
\end{DEf}

\begin{exem}\label{exe}Nous aurons besoin des descriptions explicites suivantes.

Commen\c cons par le cas o\`u le simplexe $\s=\Phi$ est maximal, i.e. $\Phi$ est associ\'e \`a une cha\^ine des r\'eseaux
$\varpi\Phi_{d-1}\subsetneq\Phi_0\subsetneq\cdots\subsetneq\Phi_{d-1}$.
On peut supposer que $\Phi_{d-1}=\OC^d$, et que sous la base canonique $\{e_0,\ldots,e_{d-1}\}$ de $\OC^d,$
$$
\Phi_i=\langle e_0,\ldots,e_i,\varpi e_{i+1},\ldots,\varpi e_{d-1}\rangle,~  \forall 0\leq i\leq
d-1.
$$
La condition impos\'ee sur $\alpha_i$ implique que $\alpha_i(e_i)$ engendre
$L_i$. Identifions $L_i$ avec $R$ en posant $\alpha_i(e_i)=1.$ Le morphisme $R$-lin\'eaire $\Pi:L_i\to L_{i+1}$ est alors donn\'e par multiplication par $c_i$, o\`u
\[
c_i=\frac{\alpha_{d-1}(e_i)}{\alpha_{d-1}(e_{i+1})},\ \ 0\leq i\leq
d-2\ \ \text{et}\ \
c_{d-1}=\frac{\varpi\alpha_{d-1}(e_{d-1})}{\alpha_{d-1}(e_0)}.
\]


Ceci nous permet d'identifier le sch\'ema formel
$\wh{\Omega}^{d-1}_{\OC,\Phi}$ au spectre formel du compl\'et\'e
$\varpi$-adique de l'anneau
$$
\OC[c_0,\ldots,c_{d-1},P_\Phi^{-1}]/(\prod c_i-\varpi),
$$
o\`u $P_\Phi=\prod P_{\mathfrak{a},i}$, $i\in\ZM/d\ZM$,
$\mathfrak{a}=(a_0,\ldots,a_{d-2})$ parcourt une classe de
repr\'esentants de $(\OC/\varpi\OC)^{d-1}$ dans $\OC^{d-1}$, et
$P_{\mathfrak{a},i}=1+a_0c_{i-1}+a_1c_{i-1}c_{i-2}+\cdots+a_{d-2}c_{i-1}\cdots
c_{1+i-d}$.

L'objet universel sur $\wh{\Omega}^{d-1}_{\OC,\Phi}$ est la donn\'ee
d'une suite de $\OC_{\wh{\Omega}^{d-1}_{\OC,\Phi}}$-modules
$$
\LC_{d-1}\To \Pi\LC_0\To\Pi\cdots\To\Pi\LC_{d-1}
$$
avec $\LC_i$ libre de base $1$ et $\Pi:\LC_i\to\LC_{i+1}$ est la
multiplication par $c_i$.

Ensuite on consid\`ere le cas o\`u $\s=[\Lambda]$ est un sommet.
Il suffit de traiter le cas o\`u $\Lambda=\OC^d$. Le sch\'ema
formel $\wh{\Omega}^{d-1}_{\OC,[\L]}$ classifie les diagrammes
commutatifs

$$
\xymatrix{
  \varpi\Lambda \ar[d]_{\alpha/\varpi} \ar@{^(->}[r] & \Lambda \ar[d]^{\alpha} \\
  L \ar[r]^{\Pi} & L   }
$$
 tels que l'application
\[
\alpha(x):\Lambda/\varpi\Lambda\To {} L\otimes_R k(x)
\]
soit injective pour tout $x\in\Spec(R).$ Ceci implique que $\alpha(u)$ est un g\'en\'erateur
de $L$, $\forall u\in\Lambda\ba\varpi\Lambda$. Le couple $(L,\alpha)$
est alors d\'etermin\'e \`a isomorphisme pr\`es par
\[
\biggl
(x_0=\frac{\alpha(e_0)}{\alpha(e_{d-1})},\ldots,x_{d-2}=\frac{\alpha(e_{d-2})}{\alpha(e_{d-1})}\biggr)\in
R^{d-1}.
\]
Nous pouvons donc identifier le sch\'ema formel
$\wh{\Omega}^{d-1}_{\OC,[\L]}$ au spectre formel du compl\'et\'e
$\varpi$-adique de l'anneau
\[
\OC[x_0,\ldots,x_{d-2},P_\Lambda^{-1}],
\]
o\`u $P_\Lambda=\prod
(a_0x_0+a_1x_1+\cdots+a_{d-2}x_{d-2}+a_{d-1})$,
$(a_0,\ldots,a_{d-1})$ parcourt une classe de repr\'esentants de
$(\OC/\varpi\OC)^d\ba \{0\}$ dans $\OC^d$.
L'objet universel sur $\wh{\Omega}^{d-1}_{\OC,[\L]}$ est la
donn\'ee d'un $\OC_{\wh{\Omega}^{d-1}_{\OC,[\L]}}$-module $\LC$
libre de base $1$, et $\Pi:\LC\to\LC$ est la multiplication par
$\varpi$.

\begin{rem}L'immersion canonique $\wh{\Omega}_{\OC,[\Lambda]}^{d-1}\into\wh{\Omega}_{\OC,\Phi}^{d-1}$
(ou
$\Omega^{d-1}_{K,[\Lambda]}\into\Omega^{d-1}_{K,\Phi}\into\PM^{d-1}_K$)
induit une identification $c_i=x_i/x_{i+1}$ pour $0\leq i\leq d-3$,
$c_{d-2}=x_{d-2}$, $c_{d-1}=\varpi/x_{0}.$
\end{rem}

Enfin, soit
$\s=\{\varpi\Lambda_{d-1}\subsetneq\Lambda_0\subsetneq\Lambda_{d-1}=\OC^{d}\}$
le sous-simplexe de $\Phi$ de type $(1,d-1),$ i.e. $\Lambda_0$ est
engendr\'e par $e_0,\varpi e_1,\ldots,\varpi e_{d-1}.$ L'objet
universel correspondant est d\'ecrit par le diagramme commutatif suivant

$$
\xymatrix{
  \varpi\Lambda \ar[d]_{\alpha/\varpi} \ar@{^(->}[r] & \Lambda_0 \ar[d]_{\alpha_0} \ar@{^(->}[r] & \Lambda \ar[d]^{\alpha} \\
  \LC \ar[r]^{\varpi/x_0} & \LC_0 \ar[r]^{x_0} & \LC   }
$$
Le sch\'ema formel $\wh{\Omega}^{d-1}_{\OC,\s}$ s'identifie au
spectre formel du compl\'et\'e $\varpi$-adique de l'anneau
\[
\OC[x_0,\ldots,x_{d-2},c_{d-1},P_{\s}^{-1}]/(x_0c_{d-1}-\varpi)
\]
o\`u
$P_{\s}=\prod(1+a_0x_0+\cdots+a_{d-2}x_{d-2})(1+a_0x_1c_{d-1}+\cdots+a_{d-3}x_{d-2}c_{d-1}+a_{d-2}c_{d-1})$,
$(a_0,\ldots,a_{d-2})$ parcourt une classe de repr\'esentants de
$(\OC/\varpi\OC)^{d-1}$ dans $\OC^{d-1}$.
\end{exem}

\ali\label{Sec:2} Les composantes irr\'eductibles de la fibre sp\'eciale g\'eom\'etrique $\o{\Omega}:=\wh{\Omega}^{d-1}_{\OC}\otimes_\OC\o{\FM}_q$ de $\widehat{\Omega}^{d-1}_{\OC}$ sont param\'etr\'ees par les sommets de $\BC\TC.$ Plus pr\'ecis\'ement, soit $s=[\L_s]$ repr\'esent\'e par un r\'eseau $\L_s,$ consid\'erons tous les simplexes maximaux contenant $s.$ Pour un tel simplexe $\s,$ il est repr\'esent\'e par une suite de r\'eseaux $\{\varpi
\Lambda_s\subsetneq\Lambda_{\s,0}\subsetneq\cdots\subsetneq\L_{\s,d-2}\subsetneq\Lambda_{s}\}.$ Notons $\o{\Omega}_s$ la vari\'et\'e projective obtenue \`a partir de $\PM(\L_s/\varpi\L_{s})$ en l'\'eclatant successivement le long du sous-sch\'ema ferm\'e $\PM(\L_s/\L_{\s,d-2})$ pour tout simplexe maximal $\s$ contenant $s,$ puis en \'eclatant le transform\'e strict de $\PM(\L_s/\L_{\s,d-3})$ pour tous ces $\s,$ puis le transform\'e strict de $\PM(\L_s/\L_{\s,d-4})$ pour tous ces $\s$ et ainsi de suite, {\em cf.} \cite[\S 4]{ito} ou \cite[(4.1.2)]{wang_DL}. On sait alors que le $\oFq$-sch\'ema $\o{\Omega}$ est localement de type fini, et $\o{\O}=\bigcup_{s\in\BC\TC_0}\o{\Omega}_s.$ Chaque $\o{\O}_s$ est une vari\'et\'e projective munie d'une action de $G_s:=\Stab_{G}(s).$ Soit $\s=\{s_0,\ldots,s_k\}$ un simplexe quelconque, notons $\o{\Omega}_\s$ la vari\'et\'e projective $\o{\Omega}_{s_0}\cap\cdots\cap\o{\Omega}_{s_k}$ munie d'une action de $\wh{G_\s}:=\Stab_{G}(\s).$ On d\'esigne $G_\s$ le fixateur de $\s$, et $G_\s^+$ le pro-$p$-radical de $G_\s$ (voir \cite{SS-crelle}). Notons $\o{\Omega}_\s^0:=\o{\Omega}_\s\backslash\bigcup_{s'\not\in \s}\o{\Omega}_{s'},$ et $j_\s:\o{\Omega}_\s^0\into\o{\Omega}_\s$ l'inclusion naturelle. En particulier,  $\o{\O}^0_s$ est la fibre sp\'eciale g\'eom\'etrique de $\wh{\Omega}^{d-1}_{\OC,s}$ (le mod\`ele de Deligne au-dessus de $s$).

Rappelons que Berkovich a d\'efini dans ce cas un morphisme de sp\'ecialisation $$\spe:\Omega^{d-1,ca}_K\To{}\o{\Omega}$$ (qui est appel\'e le {\em morphisme de r\'eduction} dans \cite[\S 1]{berk-vanishing}). L'image r\'eciproque sous ce morphisme de l'inclusion naturelle $j_s:\o{\Omega}_s^0\into\o{\Omega}_s$ s'identifie \`a:
$$
\xymatrix{
\spe^{-1}(\o{\O}^0_s) \ar @{=}[d]  \ar@{^(->}[r]^{\spe^{-1}(j_s)} &  \spe^{-1}(\o{\O}_s)  \ar@{=}[d]\\
\tau^{-1}(|s|)  \ar@{^(->}[r] & \tau^{-1}(|s|^*)
}
$$

\begin{lemme}
Soit $s$ un sommet. Quitte \`a choisir une base d'un r\'eseau qui repr\'esente $s,$ on a un isomorphisme $G_s/G_s^+\simto \GL_d(\FM_q).$ Via cet isomorphisme, $\o{\Omega}_s^0$ munie de l'action de $G_s/G_s^+$ ($G_s^+$ agit trivialement sur $\o\O_s^0$) est isomorphe \`a $\O^{d-1}_{\oFq}$ muni de l'action de $\GL_d(\FM_q),$ o\`u $\O^{d-1}_{\oFq}$ est le compl\'ementaire de tous les hyperplans $\FM_q$-rationnels dans $\PM^{d-1}_{\oFq}$ appel\'e l'espace de Drinfeld sur le corps fini $\FM_q$ ({\em cf.} \ref{Sec:1}).
\end{lemme}
\begin{preuve}
On peut supposer que $s=[\L]=[\OC^d]$ le r\'eseau standard. D'apr\`es l'exemple \ref{exe}, $\o{\O}^0_{s}=\Spec\oFq[x_0,\ldots,x_{d-2},\o{P}_\Lambda^{-1}]$ qui s'identifie donc \`a $\O^{d-1}_{\oFq}.$ Les \'el\'ements de $G^+_s=1+\varpi M_d(\OC)$ agissent bien trivialement, et les actions sont compatibles.
\end{preuve}

\subsection{Le rev\^etement mod\'er\'e $\Sig^{ca}$}\label{o}
\ali\label{new1} Rappelons tout d'abord deux descriptions modulaires de notre sch\'ema formel $\wh{\Omega}_\OC^{d-1}$ introduites par Drinfeld dans \cite{drinfeld-covering} (voir aussi \cite{boutot-carayol}). Si $R$ est une
$\OC$-alg\`ebre, nous noterons $R[\Pi]$ le quotient de l'alg\`ebre
de polyn\^omes $R[X]$ par l'id\'eal engendr\'e par $X^d-\varpi$. C'est
donc un $R$-module libre de rang $d$, engendr\'e par $1$ et un
\'el\'ement $\Pi$ (l'image de $X$) qui v\'erifie $\Pi^d=\varpi$.
L'alg\`ebre $R[\Pi]$ est munie d'une graduation \`a valeurs dans
$\ZM/d\ZM$ telle que les \'el\'ements de $R$ soient de degr\'e 0, et
$\Pi$ soit de degr\'e $1$.

On consid\`ere le foncteur $F^{Dr}$ qui associe \`a une alg\`ebre $R\in \Nilp_\OC$ l'ensemble des
classes d'isomorphie de $(\psi,\eta,T,u,r),$ o\`u
\begin{itemize}
\item $\psi$ est un $\FM_q$-homomorphisme de $\oFq$ vers $R/\varpi R.$

\item $\eta$ est un faisceau en $\OC[\Pi]$-modules plats, $\ZM/d\ZM$-gradu\'e, constructible, sur $S:=\Spec(R)$ muni de la topologie de Zariski.

\item $T$ est un faisceau en $\OC_S[\Pi]$-modules, $\ZM/d\ZM$-gradu\'e, tel que les composantes homog\`enes soient des faisceaux inversibles sur $S.$

\item $u$ est un homomorphisme $\OC[\Pi]$-lin\'eaire de degr\'e 0 de $\eta$ vers $T,$ tel que $u\otimes_\OC \OC_S:\eta\otimes_\OC\OC_S\to T$ soit surjectif.

\item $r$ est un isomorphisme $K$-lin\'eaire du faisceau constant $\underline{K}^d$ vers le faisceau $\eta_0\otimes_\OC K.$

\end{itemize}
satisfaisant les conditions suivantes:

\begin{description}
  \item[.] Soit $S_i\subset S$ le lieu d'annulation du morphisme $\Pi:T_i\to T_{i+1},$ alors la restriction $\eta_i|_{S_i}$ est un faisceau constant de fibre isomorphe \`a $\OC^d.$
  \item[.] Pour tout point $s\in S$ l'application $\eta_s/\Pi\eta_s\to (T_s/\Pi T_s)\otimes k(s)$ est injective, o\`u $k(s)$ est le corps r\'esiduel de $s$.
  \item[.] $\bigwedge^d(\eta_i)|_{S_i}=\varpi^{-i}(\bigwedge^d(\Pi^i r\underline{\OC}^d))|_{S_i}$ ($\forall i\in \ZM/d\ZM$).
\end{description}

Drinfeld d\'emontre que ce foncteur $F^{Dr}$ est pro-repr\'esentable par le
$\OC$-sch\'ema formel $\Omenr:=\wh{\Omega}^{d-1}_{\OC}\wh{\otimes}_\OC\breve{\OC}.$ Dans la suite, on notera
$(\psib,\etab,\Tb, \ub, \rb)$ l'objet universel sur
$\wh{\Omega}_{\breve{\OC}}^{d-1}$. Les composantes homog\`enes universelles $\Tb_i$ sont des fibr\'es en
droites sur $\Omenr$. Par la construction de
l'isomorphisme entre $\Omenr$ et $F^{Dr}$ ({\em cf.} \cite{boutot-carayol}), on sait que la
restriction de $\Tb_i$ \`a chaque
$\wh{\Omega}^{d-1}_{\breve{\OC},\s}:=\wh{\Omega}^{d-1}_{\OC,\s}\wh\otimes_{\OC}{\breve{\OC}}$ est en fait libre. Soit $\s=\Phi$ le simplexe maximal standard, et identifions $\Tb_i|_{\wh{\Omega}^{d-1}_{\breve{\OC},\Phi}}$ avec
$\OC_{\wh{\Omega}^{d-1}_{\breve{\OC},\Phi}}.$ Via cette identification, l'application
$\Pi:\Tb_i\to\Tb_{i+1}$ est induite par multiplication par $c_i$ ({\em cf.} l'exemple \ref{exe}).


Le foncteur $G^{Dr}$ de Drinfeld est un probl\`eme de modules des $\OC_D$-modules formels munis d'une certaine rigidification. Rappelons ci-dessous leurs d\'efinitions. Si $R$ est une $\OC$-alg\`ebre, un $\OC$-module formel $X$ est un groupe formel sur $R$ muni d'une action de $\OC$ relevant l'action naturelle sur l'espace tangent $\lie(X)$. Un $\OC_D$-module formel sur $R$ est un $\OC$-module formel muni d'une action de $\OC_D$ prolongeant l'action de $\OC.$ Un $\OC_D$-module formel $X$ est dit {\em sp\'ecial} si l'action de $\OC_d$ fait de $\lie(X)$ un $\OC_d\otimes_{\OC}R$-module localement libre de rang $1.$

La d\'efinition du foncteur $G^{Dr}$ repose sur l'existence d'un $\OC_D$-module formel sp\'ecial $\HM$ de dimension $d$ et ($\OC$-)hauteur $d^2$ sur $\oFq,$ qui est unique \`a isog\'enie pr\`es ({\em cf.} \cite{drinfeld-covering} voir aussi \cite[II Prop. 5.2]{boutot-carayol}). On consid\`ere le foncteur $G^{Dr}$ sur $\Nilp_{\OC}$ qui associe \`a $R\in\Nilp_{\OC}$ l'ensemble des classes d'isomorphie de triple $(\psi,X,\rho),$ o\`u
\begin{itemize}
\item $\psi$ est un $\FM_q$-homomorphisme de $\oFq$ vers $R/\varpi R.$

\item $X$ est un $\OC_D$-module formel sp\'ecial de hauteur $d^2$ sur $R.$

\item $\rho$ est une quasi-isog\'enie de hauteur z\'ero de $\psi^*\HM:=\HM\otimes_{\oFq,\psi}R/\varpi R$ vers $X_{R/\varpi R}.$

\end{itemize}

Un th\'eor\`eme difficile de Drinfeld nous dit qu'il existe un isomorphisme entre $G^{Dr}$ et $F^{Dr}.$ C'est-\`a-dire $G^{Dr}$ est pro-repr\'esentable par le $\OC$-sch\'ema formel $\Omenr.$

\ali\label{dd} On d\'esigne $\XG$ le $\OC_D$-module formel sp\'ecial universel de dimension
$d$ et hauteur $d^2$ sur $\Omenr$. Le morphisme $\Pi_D:\XG\to\XG$ est
une isog\'enie, son noyau $\XG[\Pi_D]$ est un sch\'ema formel en groupes fini plat de rang $q^d$ sur $\Omenr$. On
note $\Sigma:=\underline{\Isom}_{\OC_D}(\OC_D/\Pi_D\OC_D,\XG[\Pi_D]^{\rig})$ ou $\underline{\Isom}_{\OC_D}(\OC_D/\Pi_D\OC_D,\XG[\Pi_D]^{an})$ selon besoin. Par construction, $\Sig$ est un rev\^etement fini \'etale
sur $\Ome:=\Omega^{d-1}_K\wh{\otimes}_K\breve{K}$ de groupe de
Galois $(\OC_D/\Pi_D\OC_D)^\times\simeq\FM_{q^d}^\times$. On note $\Sigma^{ca}:=\Sigma\widehat{\otimes}_K\widehat{K^{ca}}\To p\Omega_K^{d-1,ca}$ la projection naturelle induite par $\XG[\Pi_D]\to \Omenr$ apr\`es une extension de scalaires. On sait que le groupe des quasi-isog\'enies de hauteur z\'ero du $\OC_D$-module formel $\HM$ vers lui-m\^eme s'identifie \`a $$G^\circ:=\Ker(\val_K\circ\det:\GL_d(K)\to K^\times),$$ o\`u $\val_K$ est la valuation normalis\'ee sur $K$ de sorte que $\val_K(\varpi)=1.$ Ceci fournit une action de $G^\circ$ sur tous les niveaux de la tour de Drinfeld. Par cons\'equent le morphisme de transition $p:\Sig^{ca}\to \Omega_K^{d-1,ca}$ est $G^\circ$-\'equivariant. Dans cet article, on s'int\'eresse \`a la cohomologie \'etale \`a support compact de $\Sig^{ca}$ au sens de Berkovich.

Par la construction pr\'ec\'edente, on a un diagramme commutatif dont toutes les fl\`eches sont $G^\circ$-\'equivariantes:
$$\xymatrix{
\Sig^{ca} \ar[d]^p \ar[rd]^\nu &\\
\O^{d-1,ca}_K\ar[r]^{\tau}& |\BC\TC|
}
$$
o\`u $\nu$ est la compos\'ee $\tau\circ p.$ Donc $\Sigma^{ca}$ admet un recouvrement par les ouverts admissibles $\nu^{-1}(|s|^*)$ o\`u $s$ parcourent les sommets de $\BC\TC.$ Soit $s$ un sommet quelconque de $\BC\TC,$ l'immersion ouverte $\nu^{-1}(|s|)\into\nu^{-1}(|s|^*)$ induit un morphisme de restriction:
\[
R\G(\nu^{-1}(|s|^*),\Lambda)\To{\res.} R\G(\nu^{-1}(|s|),\Lambda)
\]
o\`u $\L=\ZM/n$ avec $n$ un entier premier \`a $p.$

\begin{theo}\label{new2}
Le morphisme de restriction ci-dessus est en fait un isomorphisme.
\end{theo}
\begin{preuve}
La d\'emonstration repose sur un r\'esultat des cycles \'evanescents d'un faisceau mod\'er\'e sur une vari\'et\'e de r\'eduction semi-stable \'etabli par Zheng {\em cf.} \cite[Lemme 5.6]{zheng}. Supposons tout d'abord que $\G$ soit un sous-groupe discret cocompact et sans torsion de $\PGL_d(K).$ On sait alors que $\G$ agit librement sur $\Omenr$ de sorte que $\Omenr/\G$ soit propre sur $\Spf\breve{\OC}.$ On peut munir $\O^{d-1}_{\breve{K}}/\G$ d'une structure de $\breve{K}$-espace analytique telle que le quotient $\O^{d-1}_{\breve{K}}\onto \O^{d-1}_{\breve{K}}/\G$ soit un rev\^etement analytique galoisien. D'apr\`es Kurihara \cite{Kurihara} et Mustafin \cite{Mustafin}, $\O^{d-1}_{\breve{K}}/\G$ est alg\'ebrisable, i.e. il existe un sch\'ema propre $X_\G$ sur $\breve{\OC}$ de r\'eduction semi-stable tel que $\wh{\O}^{d-1}_{\breve{\OC}}/\G$ soit le compl\'et\'e de $X_\G$ le long de sa fibre sp\'eciale $X_{\G0}.$ La stratification de $X_{\G0}$ co\"incide avec le complexe simplicial $\BC\TC/\G$ {\em cf.} \cite[Thm. 2.2.6]{Kurihara}. Quitte \`a rapetisser $\G,$ on peut supposer que la projection $\pi:\Omenr\onto\Omenr/\G$ induit un isomorphisme entre $\o{\O}_s$ et $\pi(\o{\O}_s).$ Notons $Z$ (resp. $Z^0$) le sous-sch\'ema localement ferm\'e de $X_\G$ qui correspond \`a $\pi(\o{\O}_s)$ (resp. $\pi(\o{\O}^0_s)$) par l'alg\'ebrisation. D'apr\`es GAGA analytique, le rev\^etement \'etale mod\'er\'ement ramifi\'e $\Sig_{\breve{K}}/\G$ de $\O^{d-1}_{\breve{K}}/\G$ correspond \`a un rev\^etement mod\'er\'e $f:Y\onto X_{\G,\eta}$ de $X_{\G,\eta},$ o\`u $X_{\G,\eta}$ est la fibre g\'en\'erique de $X_\G.$

 D'apr\`es \cite[Corollary 3.5]{berk-vanishing}, on a des isomorphismes canoniques:
\begin{align*}
&R\G(\nu^{-1}(|s|^*),\L)=R\G(\tau^{-1}(|s|^*),p_*\L|_{\tau^{-1}(|s|^*)})\simto R\G(\o{\O}_s,R\Psi_\eta(p_*\L)|_{\o{\O}_s})\\
&R\G(\nu^{-1}(|s|),\L)=R\G(\tau^{-1}(|s|),p_*\L|_{\tau^{-1}(|s|)})\simto R\G(\o{\O}^0_s,R\Psi_\eta(p_*\L)|_{\o{\O}^0_s})
\end{align*}
qui nous donnent un diagramme commutatif
$$
\xymatrix{
R\G(\nu^{-1}(|s|^*),\Lambda)\ar[d]^\cong\ar[r]^{\res.} & R\G(\nu^{-1}(|s|),\Lambda)\ar[d]^\cong\\
R\G(\o{\O}_s,R\Psi_\eta(p_*\L)|_{\o{\O}_s})\ar[r]^{\res.} & R\G(\o{\O}^0_s,R\Psi_\eta(p_*\L)|_{\o{\O}^0_s})
}
$$
o\`u $R\Psi_\eta$ d\'esigne le foncteur des cycles \'evanescents formels d\'efini par Berkovich dans {\em loc. cit.}. Le th\'eor\`eme principal de Berkovich nous dit qu'il existe des isomorphismes canoniques $$R\Psi_\eta(p_*\L)|_{\o{\O}_s}\cong R\Psi(X_\G,f_*\L)|_Z$$ et $$R\Psi_\eta(p_*\L)|_{\o{\O}^0_s}\cong R\Psi(X_\G,f_*\L)|_{Z^0},$$ o\`u $R\Psi$ est le foncteur des cycles \'evanescents alg\'ebrique.
En vertu du \cite[Lemme 5.6]{zheng}, nous avons $$R\Psi(X_\G,f_*\L)|_Z=Rj_{Z*}R\Psi(X_\G,f_*\L)|_{Z^0},$$ o\`u $j_Z$ d\'esigne l'immersion naturelle $Z^0\into Z.$ On en d\'eduit l'\'egalit\'e suivante
\ini\begin{equation}\label{Eq::3}
 R\Psi_\eta(p_*\L)|_{\o{\O}_s}=Rj_{s,*}R\Psi_\eta(p_*\L)|_{\o{\O}^0_s},
\end{equation}
et donc un isomorphisme $$R\G(\nu^{-1}(|s|^*),\Lambda)\simto R\G(\nu^{-1}(|s|),\Lambda)$$ donn\'e par le morphisme de restriction.
\end{preuve}

\subsection{Le prolongement de $\Sigma^{ca}$ au-dessus d'un sommet}\label{S2Pa3}

Dans ce paragraphe, on prolonge le rev\^etement mod\'er\'e $\Sig^{ca}$ en un $\mu_{q^d-1}$-torseur sur $\o{\O}^0_s$ pour chaque sommet $s$ de $\BC\TC.$ Les calculs que nous effectuons ici g\'en\'eralisent ceux de Teitelbaum \cite{teitelbaum} pour $d=2.$

\ali L'espace tangent $\lie(\XG)$ du $\OC_D$-module formel sp\'ecial universel $\XG$ (voir \ref{dd}) admet une graduation par $\ZM/d\ZM$ sous l'action de $\OC_d\subset\OC_D$ en posant
$$\lie(\XG)_i=\{m\in\lie(\XG)\;|\;
\iota(a)(m)=\sigma^{-i}(a)m,\forall a\in\OC_d\}$$ o\`u
$\iota:\OC_D\to \End(\XG)$ exprime la structure de $\OC_D$-module de
$\XG$. Chaque $\lie(\XG)_i$ est un faisceau inversible sur
$\Omenr$. Consid\'erons l'objet universel $(\psib,\etab,\Tb, \ub, \rb)$ sur
$\wh{\Omega}_{\breve{\OC}}^{d-1}$ rappel\'e dans \ref{new1}. Dans
la construction de l'isomorphisme de $G^{Dr}\simto
F^{Dr}$ (voir \cite[Th\'eor\`eme 8.4]{boutot-carayol}), on identifie $\Tb_i$ \`a $\lie(\XG)_i,$ et l'action de $\Pi$ envoie $\Tb_i$ vers $\Tb_{i+1}.$ On en d\'eduit une d\'ecomposition de l'espace cotangent de $\XG[\Pi]$:
\ini\begin{equation}\label{1 dec}
\lie(\XG[\Pi])^\vee=\Tb_{0}^\vee/\Pi\Tb_{1}^\vee\oplus\Tb_{1}^\vee/\Pi\Tb_{2}^\vee\oplus\cdots\oplus\Tb_{d-1}^\vee/\Pi\Tb_{0}^\vee.
\end{equation}
Le morphisme $\iota$ induit un
homomorphisme $\o \iota: \OC_D/\Pi\simeq\FM_{q^d}\to \End(\XG[\Pi])$. Ceci nous permet d'utiliser la classification de Raynaud \cite{Raynaud} que nous rappelons ci-dessous.

Soit $M=\Hom(\FM_{q^d}^{\times},\OC_D^{\times})$ le groupe des
caract\`eres (homomorphisme de groupes) de $\FM_{q^d}^\times$ \`a valeurs  dans $\OC_D^\times$.
On prolonge chaque caract\`ere $\mu\in M$ \`a $\FM_{q^d}=\OC_D/\Pi\OC_D$
tout entier en posant $\mu(0)=0$. Un caract\`ere $\mu$ est dit {\em
fondamental} si l'application compos\'ee $\FM_{q^d}\To \mu \OC_D\To
{\can.}\FM_{q^d}$ est un homomorphisme de corps. On a donc $fd$
caract\`eres fondamentaux au total. Si on d\'esigne $\chi:\FM^{\times}_{q^d}\to\OC_d^{\times}\subset\OC_D^{\times}$ le
repr\'esentant de Teichm\"uller. Alors l'ensemble des caract\`eres fondamentaux $\{\chi_i\}_{0\leq i\leq fd-1}$ sont de la forme
$\chi_0=\chi$, $\chi_{i}=\chi_{i-1}^p$, $1\leq i\leq fd-1$. Notons
que $\bar{\iota}(\l)=\iota(\chi(\l))|_{\XG[\Pi]}\in\End(\XG[\Pi])$, pour tout $
\l\in\FM_{q^d}$.

Soit $\OC_{\XG[\Pi]}=\OC_{\Omenr}\oplus\IC$, o\`u $\IC$ est
l'id\'eal d'augmentation. L'endomorphisme $\o\iota(\l)$ sur $\XG[\Pi]$
induit un endomorphisme $[\l]$ de l'alg\`ebre de Hopf
$\OC_{\XG[\Pi]}$. Pour tout $\mu\in M$, les endomorphismes
\[
i_\mu=\frac{1}{q-1}\sum_{\l\in\FM_{q^d}^\times}\mu^{-1}(\l)[\l]
\]
de la $\OC_{\Omenr}$-alg\`ebre $\OC_{\XG[\Pi]}$ forment une famille
d'idempotents orthogonaux qui respectent $\IC$. On a alors une
d\'ecomposition
\[
\IC=\bigoplus_{\mu\in M}\IC_\mu
\]
o\`u $\IC_\mu=i_\mu(\IC)$ form\'e des \'el\'ements $x\in\IC$ tels que
 $[\l](x)=\mu(\l)x$, pour tout $\l\in\FM_{q^d}^\times.$ Notons $\IC_i:=\IC_{\chi_i}$, $\forall 0\leq i\leq fd-1$. Le $\Omenr$-sch\'ema en groupes $\XG[\Pi]$
satisfait la condition de la classification de Raynaud, i.e. chacun des faisceaux $\IC_i$ est un $\OC_{\Omenr}$-module inversible, {\em cf.} \cite[Prop. 1.2.2]{Raynaud}. Donc $\IC$ est un $\OC_{\Omenr}$-module localement
libre de rang $q^d-1$.

\begin{fact}(\cite[Thm. 1.4.1]{Raynaud}) Sous ces conditions, le sch\'ema en groupes $\XG[\Pi]$ est d\'etermin\'e par le syst\`eme $(\IC_i,\o{c}_i: \IC_{i+1}\to \IC_i^p,\o{d}_i:\IC_i^p\to\IC_{i+1})_i$ o\`u les $\o{c}_i$ et $\o{d}_i$ sont $\OC_{\Omenr}$-lin\'eaires de sorte que $\o{d}_i\circ \o{c}_i=w\id_{\IC_{i+1}}, \forall 0\leq i\leq fd-1.$ Ici $w\in\G(\Omenr,\OC_{\Omenr})$ est donn\'e par la somme de Gauss ind\'ependant du $\XG[\Pi].$
\end{fact}

Localement, on peut supposer que chaque $\IC_i$ est en fait libre engendr\'e par $X_i.$ On en d\'eduit que $\XG[\Pi]$ est donn\'e localement sur $\Omenr$ par les \'equations
$$X_i^p=\delta_iX_{i+1},i\in\ZM/fd\ZM$$ avec $\d_i$ des sections locales de $\Omenr.$
On a alors une autre description de l'espace cotangent de $\XG[\Pi]:$
\ini\begin{equation}\label{2 dec}
\lie(\XG[\Pi])^\vee=\IC/\IC^2=\IC_{0}/\IC_{fd-1}^p\oplus\IC_{1}/\IC_{0}^p\oplus\cdots\oplus\IC_{fd-1}/\IC_{fd-2}^p.
\end{equation}


\begin{lemme}
$\Tb_{i}^\vee/\Pi\Tb_{i+1}^\vee=\{x\in\IC/\IC^2\;|\; [\l](x)=\sigma^{-i}(\chi(\l))x,\forall
\l\in\FM_{q^d}\}$.
\end{lemme}
\begin{preuve}
On a $\lie(\XG[\Pi])=\Hom(\IC/\IC^2,\OC_{\Omenr}).$ Pour tout $\l\in\FM_{q^d},$ $\iota(\l)$ induit une action sur $\IC/\IC^2$ qui est, par d\'efinition, celle induite par $[\l].$
\end{preuve}

En comparant les deux d\'ecompositions \ref{1 dec} et \ref{2 dec}, on obtient
que
\[
\IC_{fi+j}=\IC_{fi+j-1}^p, 1\leq j\leq f-1, \text{ et }
\IC_{fi}/\IC_{fi-1}^p=\Tb_{d-i}^\vee/\Pi\Tb_{d-i+1}^\vee.
\]
On en d\'eduit que
\[
\IC_{fi}/\IC_{f(i-1)}^q=\Tb_{d-i}^\vee/\Pi\Tb_{d-i+1}^\vee.
\]
Il s'ensuit que $\XG[\Pi]$ est localement donn\'e par les \'equations
$X_i^q=\delta_i X_{i+1}\ \ (i\in\ZM/d\ZM)$.

\ali Pour un simplexe $\s\in\BC\TC,$ on note
$\wh{\Omega}_{\breve{\OC},\s}^{d-1}$ le produit fibr\'e de
$\wh{\Omega}_{\OC,\s}^{d-1}$ avec $\Omenr$ au-dessus de
$\wh{\Omega}_\OC^{d-1}$, et $\Omega_{\breve{K},\s}^{d-1}$ sa
fibre g\'en\'erique. Si $\s=\Phi$ est le simplexe maximal
standard, le morphisme $\Pi:\Tb_{i}\to\Tb_{i+1}$ est donn\'e
par multiplication par $c_i$ ({\em cf.} \ref{new1} et \ref{exe}). Alors, sur
$\wh{\Omega}_{\breve{\OC},\Phi}^{d-1}$
\[
\IC_{fi}/\IC_{f(i-1)}^q=\Tb_{d-i}^\vee/c_{d-i}\Tb_{d-i}^{\vee}.
\]

\begin{lemme}\label{new3}
Soit $s$ un sommet de $\BC\TC$. Notons $\Sigma_{s}$ l'espace rigide $\Sigma\times_{\Omega_{\breve{K}}^{d-1}}\Omega_{\breve{K},s}^{d-1}.$ Alors il existe une section $u\in\G(\wh{\Omega}_{\breve{\OC},s}^{d-1},\OC_{\wh{\Omega}_{\breve{\OC},s}^{d-1}}^*)$ telle que
$$
\Sigma_{s}\cong\Sp \OC_{\Omega_{\breve{K},s}^{d-1}}[X_0]/(X_0^{q^d-1}-\varpi
u).
$$
\end{lemme}
\begin{preuve}
Il suffit de traiter le cas o\`u $s$ est le sommet standard $\L=[\OC^d].$ On observe tout d'abord que le groupe de Picard de
$\wh{\Omega}_{\breve{\OC},[\Lambda]}^{d-1}$ est trivial. En effet, il est isomorphe au groupe de
Picard de sa fibre sp\'eciale, car
$\wh{\Omega}_{\breve{\OC},[\Lambda]}^{d-1}$ est $\varpi$-adique complet ({\em cf.} \cite[3.7.4]{put}). Donc
$$\Pic(\wh{\Omega}_{\breve{\OC},[\Lambda]}^{d-1})
=\Pic(\Omega^{d-1}_{\o\FM_q})=0,$$
d'apr\`es \cite[II Prop. 6.5]{hartshorne}. Donc
$\XG[\Pi]\times_{\Omenr}\wh{\Omega}_{\breve{\OC},[\Lambda]}^{d-1}$ est
donn\'e par $X_i^q=\d_i X_{i+1}$ o\`u $\d_i\in\G(\wh{\Omega}_{\breve{\OC},[\Lambda]}^{d-1},\OC_{\wh{\Omega}_{\breve{\OC},[\Lambda]}^{d-1}}).$ D'apr\`es l'exemple \ref{exe}, il existe des sections $u_i\in\G(\wh{\Omega}_{\breve{\OC},[\Lambda]}^{d-1},\OC_{\wh{\Omega}_{\breve{\OC},[\Lambda]}^{d-1}}^*)$ telles que $\d_i=c_i u_i,~\forall i.$ On en d\'eduit que
\begin{align*}
\Sigma_{[\L]}&=\Sp \OC_{\Omega_{\breve{K},[\L]}^{d-1}}[X_0]/(X_0^{q^d-1}-\d_{0}^{q^{d-1}}\d_{1}^{q^{d-2}}\cdots\d_{d-1})\\
&=\Sp \OC_{\Omega_{\breve{K},[\L]}^{d-1}}[X_0]/(X_0^{q^d-1}-c_0\cdots c_{d-1}u_0^{q^{d-1}}u_1^{q^{d-2}}\cdots u^q_{d-2} u_{d-1}c_{0}^{q^{d-1}-1}c_{1}^{q^{d-2}-1}\cdots c^{q-1}_{d-2})\\
&=\Sp \OC_{\Omega_{\breve{K},[\L]}^{d-1}}[X_0]/(X_0^{q^d-1}-\varpi u),
\end{align*}
o\`u $u:=u_0^{q^{d-1}}u_1^{q^{d-2}}\cdots u_{d-1}c_{0}^{q^{d-1}-1}c_{1}^{q^{d-2}-1}\cdots c^{q-1}_{d-2}\in\G(\wh{\Omega}_{\breve{\OC},[\L]}^{d-1},\OC_{\wh{\Omega}_{\breve{\OC},[\L]}^{d-1}}^*),$ car $c_0\cdots c_{d-1}=\varpi$ et $c_0,\ldots,c_{d-2}$ appartiennent \`a $\G(\wh{\Omega}_{\breve{\OC},[\L]}^{d-1},\OC_{\wh{\Omega}_{\breve{\OC},[\L]}^{d-1}}^*),$ {\em cf.} l'exemple \ref{exe}. D'o\`u l'\'enonc\'e du lemme.
\end{preuve}

\ali\label{Sec:4} Posons $\breve{K}^t=\breve{K}[\varpi_t]/(\varpi_t^{q^d-1}-\varpi)$ une
extension mod\'er\'ement ramifi\'ee de degr\'e $q^d-1$ de
$\breve{K}$, $\breve{\OC}^t$ l'anneau des entiers de $\breve{K}^t$.
Apr\`es l'extension des scalaires \`a $\breve{K}^t$, l'espace rigide $\Sigma_{s,\breve{K}^t}:=\Sigma_{s}\otimes_{\breve{K}}\breve{K}^t$
est donn\'e par
$$\Sp\OC_{\Omega^{d-1}_{\breve{K}^t,s}}[X_0']/(X_0'^{q^d-1}-u),$$
o\`u $X_0'=X_0/\varpi_t$ et
$\Omega^{d-1}_{\breve{K}^t,s}=\Omega^{d-1}_{\breve{K},s}\otimes_{\breve{K}}\breve{K}^t
$. Consid\'erons la normalisation de
$\wh{\Omega}_{\breve{\OC},s}^{d-1}$ dans
$\Sigma_{s,\breve{K}^t}$ que l'on notera
$\wh{\Sigma}^0_{s}.$ D'apr\`es \cite[Exp. I (9.10), (10.2)]{SGA1}, $\wh{\Sigma}^0_{s}=\Spf\OC_{\wh{\Omega}^{d-1}_{\breve{\OC}^t,s}}[X_0']/(X_0'^{q^d-1}-u)$,
avec $u\in\G(\wh{\Omega}^{d-1}_{\breve{\OC}^t,s},\OC_{\wh{\Omega}^{d-1}_{\breve{\OC}^t,s}}^*).$ Le groupe de Galois $\Gal(\breve{K}^t/\breve{K})$ est isomorphe canoniquement \`a $\mu_{q^{d}-1}$ via $g\in\Gal(\breve{K}^t/\breve{K})\mapsto g(\varpi_t)/\varpi_t\in\mu_{q^{d}-1}.$ On en d\'eduit qu'un \'el\'ement $\z\in\mu_{q^{d}-1}$ agit sur l'anneau de sch\'ema affine formel $\wh{\Sigma}^0_{s}$ en envoyant $X'_0$ vers $\z^{-1}X'_0.$

Notons que la fibre speciale
$\o{\Sigma}_{s}^0$ de
$\wh{\Sigma}^0_{s}$ est un
$\mu_{q^d-1}$-torseur $G_s/G_s^+$-invariant au-dessus de
$\o{\O}^0_s,$ car $\o u:=u \pmod{\varpi_t}$ est
une unit\'e dans $\G(\o{\O}^0_s,\OC_{\o{\O}^0_s}).$ On a alors un diagramme commutatif:

$$\xymatrix{
\Sigma_{s,\breve{K}^t}\ar[d]_p \ar@{^(->}[r] &\wh{\Sigma}^0_{s} \ar@{<-^)}[r] \ar[d]_{\wh{p}_s} & \o{\Sigma}_{s}^0 \ar[d]_{\o{p}_s}\\
\Omega^{d-1}_{\breve{K}^t,s} \ar@{^(->}[r] &\wh{\Omega}_{\breve{\OC},s}^{d-1} \ar@{<-^)}[r] & \o{\O}^0_{s}
}$$
o\`u $\wh{p}$ et $\o{p}$ d\'esignent les projections naturelles.

\begin{coro}\label{new4}
On a un isomorphisme
$$
R\G(\nu^{-1}(|s|),\L)\cong R\G(\o{\Sigma}_{s}^0,\L).
$$
\end{coro}
\begin{preuve}
D'apr\`es Berkovich, on a
$$
R\G(\nu^{-1}(|s|),\L)=R\G(\tau^{-1}(|s|),p_*\L)=R\G(\o{\O}^0_s,R\Psi_\eta(p_*\L)|_{\o{\O}^0_s}).
$$
Notons que $\wh{\Sigma}^0_{s}$ est un mod\`ele lisse de $\Sigma_{s,\breve{K}^t}.$ D'apr\`es \cite[Exp. I 2.4]{SGA7-1},
\ini\begin{align}\label{Eq::4}
R^0\Psi_\eta(p_*\Lambda)|_{\o{\Omega}_s^0} &=\o{p}_{s*}\Lambda,\\
R^n\Psi_\eta(p_*\Lambda)|_{\o{\Omega}_s^0} &=0, \forall n\geq 1. \notag
\end{align}
Il s'ensuit que
$$
R\G(\nu^{-1}(|s|),\L)=R\G(\o{\O}^0_s,\o{p}_{s*}\Lambda)=R\G(\o{\Sigma}_{s}^0,\L).
$$
\end{preuve}

\subsection{Un calcul du torseur $\o{\Sigma}_{s}^0$}
Dans le paragraphe pr\'ec\'edent, on a reli\'e la cohomologie du tube $\nu^{-1}(|s|)$ \`a la cohomologie de $\o{\Sigma}_{s}^0.$ Rappelons que $\o{\Sigma}_{s}^0$ est un $\mu_{q^d-1}$-torseur $G_s/G_s^+$-\'equivariant sur $\o{\O}^0_s.$ Dans ce paragraphe, notre but est de calculer sa classe dans $H^1_{et}(\o{\O}^0_s,\mu_{q^d-1}).$ Supposons d\'esormais que $s$ soit le sommet standard $\L=[\OC^d].$ On se ram\`ene donc au cas o\`u $\o{\O}^0_s=\O^{d-1}_{\oFq}$ sur lequel $\o{\Sig}^0_{[\L]}$ est un $\mu_{q^d-1}$-torseur $\GL_d(\FM_q)$-\'equivariant.

\ali\label{p} Notons $\HC$ l'ensemble des hyperplans $\FM_q$-rationnels de $\PM^{d-1}_{\oFq}$, nous avons alors
$$\Omega_{\oFq}^{d-1}=\PM^{d-1}_{\oFq}\ba \bigcup_{Y\in\HC}Y.$$
Notons $i$ (resp. $j$) l'inclusion naturelle de $D:=\bigcup_{Y\in\HC}Y$ (resp. $\Omega_{\oFq}^{d-1}$) dans $\PM^{d-1}_{\oFq}.$ Pour $I$ un sous-ensemble de $\HC$, on notera $Y_I=\bigcap_{Y\in I}Y$, et $i_{Y_I}$ l'inclusion de $Y_I$ dans $\PM^{d-1}_{\oFq}.$ La suite exacte de cohomologie relative associ\'ee aux inclusions:
$$
\xymatrix{
  \Omega_{\oFq}^{d-1} \ar@{^(->}[r]^{j} & \PM^{d-1}_{\oFq} \ar@{<-^)}^{i}[r] & \bigcup_{Y\in\HC}Y
  }
$$
nous fournit une suite exacte

\ini\begin{equation}\label{omega}
0\To {} H^1(\Omega_{\oFq}^{d-1},\mu_n)\To \partial
H^2_{D}(\PM^{d-1}_{\oFq},\mu_n)\To {}
H^2(\PM^{d-1}_{\oFq},\mu_n)
\end{equation}
o\`u $n$ est un entier premier \`a $p.$

\begin{lemme} On a un isomorphisme canonique
\[
H^2_{D}(\PM^{d-1}_{\oFq},\mu_n)=\bigoplus_{Y\in\HC}H_Y^2(\PM^{d-1}_{\oFq},\mu_n)=\bigoplus_{Y\in\HC}\ZM/n.
\]
\end{lemme}
\begin{preuve}
On prend une r\'esolution injective $\ZM/n\to\IC^\cdot$ du faisceau constant $\ZM/n.$ Pour chaque $q,$ nous avons une r\'esolution simpliciale de $i_*i^!\IC^q:$
$$
0\to\cdots\to\bigoplus_{I\subset\HC,|I|=r}i_{Y_I,*}i_{Y_I}^!(\IC^q)\to\cdots\to\bigoplus_{Y\in\HC}i_{Y,*}i_{Y}^!(\IC^q)\to i_*i^!(\IC^q)\to0.
$$
Les suites spectrales associ\'ees \`a deux filtrations du double complexe $K^{pq}:=\bigoplus_{I\subset\HC, |I|=-p}i_{Y_I,*}i^!_{Y_I}\IC^q~(p\leq -1,q\geq0)$ nous fournit une suite spectrale:
$$
E_1^{pq}=\bigoplus_{I\subset\HC, |I|=-p}i_{Y_I,*}R^qi^!_{Y_I}\ZM/n\ddonc i_*R^{p+q}i^!\ZM/n.
$$

Pour chaque hyperplan rationnel $Y,$ $(Y,\PM_{\oFq}^{d-1})$ est un couple lisse (\cite[Exp. XVI]{SGA4-3}) de codimension $1$, et $(Y_I,\PM_{\oFq}^{d-1})$ est un couple lisse de codimension $>1$ si $|I|\geq2$. D'apr\`es la puret\'e (voir {\em loc. cit.}), on sait que $$R^0i^!(\ZM/n)=R^1i^!(\ZM/n)=0,$$ et $$i_*R^2i^!(\ZM/n)=\bigoplus_{Y\in\HC}i_{Y*}R^2i^!_Y(\ZM/n)=\bigoplus_{Y\in\HC}i_{Y*}(\ZM/n) _{Y}(-1).$$ On d\'eduit la premi\`ere \'egalit\'e par la suite spectrale
$$
E^{pq}_2=H^p(D,R^qi^!\mu_n)\Longrightarrow H^{p+q}_D(\PM_{\oFq}^{d-1},\mu_n).
$$

Chaque $Y$ est un diviseur irr\'eductible, et la deuxi\`eme \'egalit\'e est donn\'ee par la classe fondamentale de $Y$:
\[
H_Y^2(\PM^{d-1}_{\oFq},\mu_n)=H^0(Y,R^2i_Y^!\mu_n)=H^0(Y,(\ZM/n)_Y)=\ZM/n.
\]
\end{preuve}

\ali Posons $n=q^d-1$ et consid\'erons la suite exacte \ref{omega}. Re\'ecrivons-la sous la forme suivante:
\[
0\To{}
H_{et}^1(\Omega_{\oFq}^{d-1},\mu_{q^d-1})\To\partial\bigoplus_{Y\in\HC}\ZM/(q^d-1)\To
\sum\ZM/(q^d-1).
\]
Notre but est de calculer la classe du $\mu_{q^d-1}$-torseur $\o{\Sigma}^0_{[\Lambda]}$ dans $H_{et}^1(\Omega_{\oFq}^{d-1},\mu_{q^d-1}).$ Rappelons tout d'abord la d\'efinition de l'application $\partial$. Comme
$\Pic(\Omega_{\oFq}^{d-1})=0$, par la suite exacte de Kummer, un $\mu_{q^d-1}$-torseur $Z$
peut \^etre \'ecrit sous la forme $\OC_{\Omega_{\oFq}^{d-1}}[T]/(T^{q^d-1}-f)$,
avec $f\in\G(\Omega_{\oFq}^{d-1},\OC_{\Omega_{\oFq}^{d-1}}^*).$ Alors $$\partial(Z)(Y)\equiv\ord_Y f \pmod{q^d-1}, ~\forall Y\in\HC.$$


\begin{lemme}\label{1}
Soit $Z$ un $\mu_{q^d-1}$-torseur $\GL_d(\FM_q)$-invariant sur
$\Omega_{\oFq}^{d-1}$, alors $\partial(Z)(Y)=\partial(Z)(g\cdot Y), \forall
g\in\GL_d(\FM_q)$, et $\partial(Z)(Y)\equiv0\pmod{q-1}$.
\end{lemme}
\begin{preuve}
La premi\`ere assertion d\'ecoule d'invariance sous $g\in
\GL_d(\FM_q)$. Pour la deuxi\`eme, on observe que $\GL_d(\FM_q)$
agit transitivement sur $\HC$ et le cardinal de $\HC$ est
$1+q+q^2+\cdots+q^{d-1}$. Notons que $\partial(Z)$ est contenu dans le
noyau de $\sum:\bigoplus_{Y\in\HC}\ZM/(q^d-1)\to\ZM/(q^d-1) $, on a donc
$$(1+q+\cdots+q^{d-1})\cdot\partial(Z)(Y)\equiv0\pmod{q^d-1}.$$ Donc $\partial(Z)(Y)\equiv0\pmod{q-1}.$
\end{preuve}

\begin{theo}\label{thm1}
Pour tout $Y\in\HC,$ on a $\partial(\o{\Sigma}^0_{[\Lambda]})(Y)\equiv q-1\pmod{q^d-1}.$
\end{theo}
\begin{preuve}
Commen\c cons par un lemme g\'eom\'etrique.
\begin{lem}
Soit $\s=\{\varpi\Lambda\subsetneq\Lambda_0\subsetneq\Lambda\}$ le
simplexe que l'on a \'etudi\'e dans l'exemple \ref{exe}, i.e.
$\Lambda_0$ correspond \`a l'hyperplan $\o x_0=0$ de
$\PM(\Lambda/\varpi\Lambda),$ alors
$\Pic(\wh{\Omega}_{\breve{\OC},\s}^{d-1})=0.$
\end{lem}

\begin{preuve}
Il suffit de montrer que le groupe de Picard de la fibre sp\'eciale $X$
de $\wh{\Omega}_{\breve{\OC},\s}^{d-1}$ est triviale. D'apr\`es l'exemple \ref{exe}, $X$ est une r\'eunion des deux composantes
irr\'eductibles
$$C=\Spec{\o\FM_q[x_0,\ldots,x_{d-2}][\prod(1+a_0x_0+\cdots+a_{d-2}x_{d-2})^{-1}]}$$
et
$$D=\Spec{\o\FM_q[x_1,\ldots,x_{d-2},c_{d-1}][P^{-1}]}$$
o\`u $P=\prod(1+a_1x_1+\cdots+a_{d-2}x_{d-2})(1+a_0x_1c_{d-1}+\cdots+a_{d-3}x_{d-2}c_{d-1}+a_{d-2}c_{d-1}),$
avec l'intersection $$E:=C\cap
D=\Spec{\o\FM_q[x_1,\ldots,x_{d-2}][\prod(1+a_1x_1+\cdots+a_{d-2}x_{d-2})^{-1}]}.$$
Nous noterons $i_C,$ $i_D,$ et $i_E$ les inclusions canoniques dans $X$. D'apr\`es \cite[II Prop. 6.5]{hartshorne}, $\Pic(C)=\Pic(D)=\Pic(E)=0.$

On montre qu'il existe une suite exacte des faisceaux sur $X$:
\ini\begin{equation}\label{picard}
1\To{}\OC_X^*\To\alpha {i_C} _*\OC_C ^*\times{i_D}
_*\OC_D^*\To\beta{i_E} _*\OC_E^*\To{}1,
\end{equation}
o\`u $\alpha$ est donn\'e par $f\mapsto(f|_C,f|_D)$, $\beta$ est
donn\'e par $(f,g)\mapsto f|_E\cdot g|_E^{-1}$. L'exactitude est
v\'erifi\'ee en regardant la fibre en chaque point. En effet, soit $x$ un
point de $X$ contenu dans $C\ba D$ (resp. $D\ba C$), le complexe
\ref{picard} se r\'eduit \`a l'isomorphisme
$\OC_{X,x}^*\simeq\OC_{C,x}^*$ (resp.
$\OC_{X,x}^*\simeq\OC_{D,x}^*$). Il nous reste le cas o\`u $x$ est
un point de $E$. Comme cette question est locale, on peut supposer
que
\begin{align*}
X&=\Spec{\o\FM_q[x_0,\ldots,x_{d-1}]/x_0x_{d-1}} &
C&=\Spec{\o\FM_q[x_0,\ldots,x_{d-1}]/x_0}\\
D&=\Spec{\o\FM_q[x_0,\ldots,x_{d-1}]/x_{d-1}} &
E&=\Spec{\o\FM_q[x_0,\ldots,x_{d-1}]/(x_0,x_{d-1})},
\end{align*}
et $x=\mathfrak{p}$ est un point de $E$. Notons que l'on a une suite
exacte des faisceaux coh\'erents sur $X$
$$
0\To{}\OC_X\To{\alpha'} {i_C} _*\OC_C \times{i_D}
_*\OC_D\To{\beta'}{i_E} _*\OC_E\To{}0
$$
o\`u $\alpha'$ est donn\'e par $f\mapsto(f|_C,f|_D)$, $\beta'$ est
donn\'e par $(f,g)\mapsto f|_E- g|_E$. On la v\'erifie en regardant
les fibres en tous les points ferm\'es. Donc on a une suite exacte
$$
0\To{}\OC_{X,\mathfrak{p}}\To{\alpha'_\mathfrak{p}}\OC_{C
,\mathfrak{p}}\times\OC_{D,\mathfrak{p}}\To{\beta'_{\mathfrak{p}}}\OC_{E,\mathfrak{p}}\To{}0.
$$
On en d\'eduit que la suite
$$
1\To{}\OC_{X,\mathfrak{p}}^*\To{\alpha_{\mathfrak{p}}}\OC_{C,\mathfrak{p}}
^*\times\OC_{D,\mathfrak{p}}^*\To{\beta_{\mathfrak{p}}}\OC_{E,\mathfrak{p}}^*
$$
est exacte. Notons que l'application surjective $$\o\FM_q[x_0,\ldots,x_{d-1}]/x_{d-1}\onto \o\FM_q[x_0,\ldots,x_{d-1}]/(x_0,x_{d-1})$$ est scind\'ee, et donc induit une surjection $\OC_{D,\mathfrak{p}}^*\onto \OC_{E,\mathfrak{p}}^*.$ C'est-\`a-dire $\beta_\mathfrak{p}$ est surjective.

On associe la suite exacte
longue de cohomologie $H^\bullet(X,-)$ au complexe \ref{picard}, et on obtient une suite
exacte
$$
H^0(C,\OC_C^*)\times H^0(D,\OC_D^*)\To{H^0(\beta)}
H^0(E,\OC^*_E)\To{} \Pic(X)\To{} \Pic(C)\times \Pic(D).
$$
L'immersion $E\into D$ est induite par l'\'epimorphisme
\begin{align*}
\varphi:\o\FM_q[x_1,\ldots,x_{d-2},c_{d-1}][P^{-1}]&\to\o\FM_q[x_1,\ldots,x_{d-2}][\prod(1+a_1x_1+\cdots+a_{d-2}x_{d-2})^{-1}]\\
x_i&\mapsto x_i\\
c_{d-1}&\mapsto 0.
\end{align*}
Notons que $\varphi$ est scind\'e, donc il induit un \'epimorphisme $H^0(D,\OC^*_D)\onto H^0(E,\OC^*_E).$ Alors $H^0(\beta)$ est surjectif, donc $\Pic(X)$
est trivial.
\end{preuve}

Revenons \`a la preuve du th\'eor\`eme. Rappelons que $$\Sigma_{[\L]}\cong\Sp \OC_{\Omega_{\breve{K},[\L]}^{d-1}}[T]/(T^{q^d-1}-\varpi u)$$ {\em cf.} le lemme \ref{new3}. Gr\^ace au lemme \ref{1}, il s'agit de calculer $\ord_H\o u$, pour un hyperplan $\FM_q$-rationnel
$H$. On peut supposer que $H$ est l'hyperplan donn\'e par $\o
x_0=0$. D'apr\`es le lemme pr\'ec\'edent, $\XG[\Pi]\times_{\Omenr}\wh{\Omega}_{\breve{\OC},\s}^{d-1}$ est donn\'e sur $\wh{\Omega}_{\breve{\OC},\s}^{d-1}$ par  les equations $X_{i}^q=\d_i X_{i+1}$ o\`u $\d_i\in\G(\wh{\Omega}_{\breve{\OC},\s}^{d-1},\OC_{\wh{\Omega}_{\breve{\OC},\s}^{d-1}}).$ En tenant compte de l'exemple \ref{exe},
il est
donn\'e par
\[
\Spf\OC_{\wh{\Omega}_{\breve{\OC},\s}^{d-1}}[X_0,X_1]/(X_0^{q^{d-1}}-\delta_1
X_1,X_1^q-\delta_0X_0)\] o\`u $\delta_1=u_1x_0$,
$\delta_0=u_0\varpi/x_0$, avec
$u_0,u_1\in\G(\wh{\Omega}_{\breve{\OC},\s}^{d-1},\OC_{\wh{\Omega}_{\breve{\OC},\s}^{d-1}}^*).$ Ceci implique que
$\Sigma_{\s}:=\Sigma\times_{\Omega_{\breve{K}}^{d-1}}\Omega_{\breve{K},\s}^{d-1}$
est l'espace rigide d\'efini par
$$ \Sp\OC_{\Omega_{\breve{K},\s}^{d-1}}[X_0]/(X_0^{q^d-1}-\varpi
u_1^qu_0x_0^{q-1}).$$

Comme $\Sigma_{\s}|_{\tau^{-1}(\Lambda)}$ et $\Sigma_{[\Lambda]}$ donnent la m\^eme classe de
torseur dans $H^1(\Omega_{\breve{K},[\Lambda]}^{d-1},\mu_{q^d-1})$, il
existe alors une fonction
$f\in\G(\Omega_{\breve{K},[\Lambda]}^{d-1},\OC_{\Omega_{\breve{K},[\Lambda]}^{d-1}}^*)$ telle que $\varpi
u=\varpi u_1^qu_0x_0^{q-1}f^{q^d-1}$. Notons que $f$ appartient \`a $\G(\wh{\Omega}_{\breve{\OC},[\Lambda]}^{d-1},\OC_{\wh{\Omega}_{\breve{\OC},[\Lambda]}^{d-1}}^*),$ en fait $u,u_0,u_1,x_0$ sont tous
dans $\G(\wh{\Omega}_{\breve{\OC},[\Lambda]}^{d-1},\OC_{\wh{\Omega}_{\breve{\OC},[\Lambda]}^{d-1}}^*).$ Donc on a
$$
\ord_{\o x_0}\o u\equiv \ord_{\o x_0}\o{u}_0+\ord_{\o x_0}\o{u}^q_1+\ord_{\o x_0}\o{x}^{q-1}_0\pmod{q^d-1}.
$$
Notons que
$u_0,u_1\in\G(\wh{\Omega}^{d-1}_{\breve{\OC},\s},\OC_{\wh{\Omega}^{d-1}_{\breve{\OC},\s}}^*)$, l'ordre de $\o u_i$ en $\o x_0$ est nulle. On en d\'eduit que
$$\ord_{\o x_0}\o u\equiv q-1\pmod{q^d-1}.$$
\end{preuve}

\subsection{Le lien avec les vari\'et\'es de Deligne-Lusztig}

Consid\'erons $\GL_d$ le groupe lin\'eaire sur $\oFq$ muni un morphisme de Frobenius $(a_{i,j})\mapsto (a_{i,j}^q).$ Dans ce paragraphe, on rappelle certains aspects de la th\'eorie de vari\'et\'es de Deligne-Lusztig \cite{deligne-lusztig} associ\'ees \`a $\GL_d$ et l'\'el\'ement de Coxeter $w:=(1,\ldots,d)\in\SG_d,$ et on \'etudie le lien avec l'espace de Drinfeld $p$-adique.

\ali\label{Sec:1}  Les vari\'et\'es de Deligne-Lusztig sont des vari\'et\'es sur un corps fini $\FM_q,$ qui jouent un r\^ole important dans l'\'etude de la th\'eorie de repr\'esentation de groupes finis r\'eductifs. Ici on s'int\'eresse au cas associ\'e \`a $\GL_d$ et l'\'el\'ement de Coxeter $w:=(1,\ldots,d)\in\SG_d,$ ceci est trait\'e dans \cite[(2.2)]{deligne-lusztig}. On a vu dans \ref{p} que la sous-vari\'et\'e ouverte $\O^{d-1}_{\FM_q}$ de $\PM^{d-1}_{\FM_q}$ est d\'efinie comme le compl\'ementaire de tous les hyperplans $\FM_q$-rationnels, donc elle est d\'efinie par la non-nullit\'e du d\'eterminant $\det((X_i^{q^j})_{0\leq i,j\leq d-1}),$ car cette fonction s'identifie (\`a une constante non-nulle pr\`es) au produit de formes lin\'eaires \`a coefficients dans $\FM_q,$ i.e.
$$
\det((X_i^{q^j})_{0\leq i,j\leq d-1})=c\cdot \prod_{[a_0:\cdots:a_{d-1}]\in\PM^{d-1}(\FM_q)}(a_0X_0+\cdots+a_{d-1}X_{d-1}),
$$
o\`u $c\in\FM_q^\times$ d\'epend du choix des rel\`evements \`a $\FM_q^d\ba\{0\}$ des \'el\'ements de $\PM^{d-1}(\FM_q).$ \'Evidemment, le groupe fini $\GL_d(\FM_q)$ agit sur $\O^{d-1}_{\FM_q}$ par translation lin\'eaire sur les coordonn\'ees.

Par la construction de Deligne et Lusztig, $\O^{d-1}_{\FM_q}$ admet un rev\^etement fini \'etale $\DL^{d-1}$ de groupe de Galois $\FM_{q^d}^\times.$ On peut identifier $\DL^{d-1}$ avec la sous-vari\'et\'e ferm\'ee de l'espace affine $\AM^d_{\FM_q}=\Spec {\FM_q[X_0,\ldots,X_{d-1}]}$ d\'efinie par l'\'equation
\ini\begin{equation}\label{Eq:2}
\det((X_i^{q^j})_{0\leq i,j\leq d-1})^{q-1}=(-1)^{d-1}.
\end{equation}
L'action de $\GL_d(\FM_q)$ sur $\O^{d-1}_{\FM_q}$ se l\`eve naturellement sur $\DL^{d-1},$ et le groupe $\FM_{q^d}^\times$ agit sur $\DL^{d-1}$ par multiplication sur les coordonn\'ees $X_i\mapsto\z X_i,~\forall \z\in\FM_{q^d}^\times.$

Posons $\DL^{d-1}_{\oFq}:=\DL^{d-1}\otimes\oFq$ le $\mu_{q^d-1}$-torseur $\GL_d(\FM_q)$-invariant au-dessus de $\Omega_{\oFq}^{d-1},$ et notons
$\pi:\DL^{d-1}_{\oFq}\to \Omega_{\oFq}^{d-1}$ la projection canonique $\GL_d(\FM_q)$-\'equivariante. Notons $x_i=X_i/X_{d-1}$, $0\leq i\leq d-2$ les coordonn\'es affines associ\'ees, alors $\O^{d-1}_{\oFq}$ s'identifie \`a $\AM_{\oFq}^{d-1}$ priv\'e les hyperplans
rationnels $a_0x_0+\cdots+a_{d-2}x_{d-2}+a_{d-1}$, o\`u
$(a_0,\ldots,a_{d-1})$ parcourt $\FM_q^{d}  \ba \{0\}$. De plus, nous avons
une expression explicite du torseur $\DL^{d-1}_{\oFq}$ en posant $T:=1/X_{d-1}$:
\begin{align*}
\DL^{d-1}_{\oFq}&=\Spec
{\OC_{\Omega_{\oFq}^{d-1}}[T]/\big(T^{q^d-1}-(-1)^{d-1}(\prod_{[a_0:\cdots:a_{d-1}]\in\PM^{d-1}(\FM_q)}(a_0x_0+\cdots+a_{d-2}x_{d-2}+a_{d-1}))^{q-1}\big)}\\
&=\Spec
{\OC_{\Omega_{\oFq}^{d-1}}[T]/(T^{q^d-1}-(-1)^d\prod_{(a_0,\ldots,a_{d-1})\in\FM_q^d\ba\{0\}}(a_0x_0+\cdots+a_{d-2}x_{d-2}+a_{d-1}))},
\end{align*}
car
$\prod_{(a_0,\ldots,a_{d-1})\in\FM_q^d\ba\{0\}}(a_0x_0+\cdots+a_{d-2}x_{d-2}+a_{d-1})=(-1)\cdot(\prod_{[a_0:\cdots:a_{d-1}]\in\PM^{d-1}(\FM_q)}(a_0x_0+\cdots+a_{d-2}x_{d-2}+a_{d-1}))^{q-1}.$ On en d\'eduit que $\z\in\FM_{q^d}^\times$ agit sur $\DL^{d-1}$ par la formule $T\mapsto \z^{-1}T.$

\begin{prop}{\label{DL}}
Rappelons que $\HC$ est l'ensemble des hyperplans $\FM_q$-rationnels de $\PM^{d-1}_{\oFq},$ {\em cf.} \ref{p}. Alors pour tout $Y\in\HC$, on a $\partial(\DL^{d-1}_{\oFq})(Y)\equiv q-1\pmod{q^d-1}$.
\end{prop}
\begin{preuve}
Le torseur $\DL^{d-1}_{\oFq}$ est $\GL_d(\FM_q)$-invariant. D'apr\`es le
lemme \ref{1}, il suffit de savoir la valeur de
$\partial(\DL^{d-1}_{\oFq})$ en l'hyperplan $X_0=0$ qui correspond \`a
l'hyperplan $x_0=0$ de $\AM_{\FM_q}^{d-1}$. Alors,
\begin{align*}
\partial(\DL^{d-1}_{\oFq})(X_0=0)&=\ord_{x_0}((-1)^d\cdot\prod_{(a_0,\ldots,a_{d-1})\in\FM_q^d\ba\{0\}}(a_0x_0+\cdots+a_{d-2}x_{d-2}+a_{d-1}))\\
&\equiv
q-1\pmod{q^d-1}.
\end{align*}
\end{preuve}

\begin{theo}\label{prop2}
Soit $s$ un sommet de $\BC\TC.$ Quitte \`a choisir une base d'un r\'eseau qui repr\'esente $s,$ on a un isomorphisme $G_s/G_s^+\cong \GL_d(\FM_q)$-\'equivariant
\ini\begin{equation}
R\G(\nu^{-1}(|s|),\Lambda)\simto R\G(\DL^{d-1}_{\oFq},\Lambda).
\end{equation}
\end{theo}
\begin{preuve}
D'apr\`es le th\'eor\`eme \ref{thm1} et la proposition \ref{DL}, on a un isomorphisme $G_s/G_s^+\cong \GL_d(\FM_q)$-\'equivariant $\o{\Sig}^0_s\cong \DL^{d-1}_{\oFq}$ et un diagramme commutatif:

$$
\xymatrix{
\o{\Sig}^0_s\ar[d]^{\o{p}_s}\ar[r]^\cong & \DL^{d-1}_{\oFq}\ar[d]^{\pi}\\
\o{\O}^0_s \ar[r]^\cong & \O^{d-1}_{\oFq}
}
$$
En vertu du corollaire \ref{new4}, on obtient un isomorphisme $$R\G(\nu^{-1}(|s|),\Lambda)\simto R\G(\DL^{d-1}_{\oFq},\Lambda).$$
\end{preuve}

\begin{coro}\label{C1C1}
On a un isomorphisme $G_s/G_s^+\cong \GL_d(\FM_q)$-\'equivariant:
$$
R\G_c(\nu^{-1}(|s|^*),\Lambda)\simto R\G_c(\DL^{d-1}_{\oFq},\Lambda).
$$
\end{coro}
\begin{preuve}
Ceci d\'ecoule du th\'eor\`eme \ref{new2} et le th\'eor\`eme pr\'ec\'edent, en vertu de la dualit\'e de Poincar\'e analytique \cite[(7.3)]{berk-ihes} et alg\'ebrique \cite[Exp. XVIII]{SGA4-3}.
\end{preuve}

\ali Dans \cite[\S 9]{deligne-lusztig}, Deligne et Lusztig ont introduit une compactification $\o \O^{d-1}_{\FM_q}$ de $\O^{d-1}_{\FM_q}$ dont le compl\'ementaire est un diviseur \`a croisements normaux. Soit $\D=\{\a_1,\ldots,\a_{d-1}\}$ l'ensemble de racines simples tel que $w=s_{\a_1}\cdots s_{\a_{d-1}}\in\SG_d,$ suivant eux, on d\'efinit
$$
\o \O^{d-1}_{\FM_q}:=\mathop{\coprod_{x=x_1\cdots x_{d-1}\in W}}_{x_i\in\{1,s_{\alpha_i}\}}X(x),
$$
o\`u $X(x)$ est la vari\'et\'e de Deligne-Lusztig associ\'ee \`a $\GL_d$ et l'\'el\'ement $x\in\SG_d,$ en particulier, $X(w)=\O^{d-1}_{\FM_q}.$ D'ailleurs, la vari\'et\'e $\o \O^{d-1}_{\FM_q}$ peut s'obtenir par une suite d'\'eclatements successifs de l'espace projectif $\PM^{d-1}_{\FM_q}$ comme dans \ref{Sec:2}, {\em cf.} \cite[Lemme 4.1.2]{wang_DL}.

\begin{fait}\label{new6}
Soit $s$ un sommet de $\BC\TC.$ Quitte \`a choisir une base d'un r\'eseau qui repr\'esente $s,$ on a un diagramme commutatif $G_s/G_s^+\cong \GL_d(\FM_q)$-\'equivariant:
$$
\xymatrix{
\o{\O}^0_s\ar[r]^\cong\ar[d] & \O^{d-1}_{\oFq}\ar[d]\\
\o{\O}_s\ar[r]^\cong & \o{\O}^{d-1}_{\oFq}.
}
$$
\end{fait}
\begin{preuve}
Ceci d\'ecoule de \cite[Lemme 4.1.2]{wang_DL}.
\end{preuve}

\begin{rem}
Soient $I$ un sous-ensemble de $\D,$ $P_I=BW_I B$ le \para standard associ\'e \`a $I$ et $U_I$ le radical unipotent de $P_I.$ D\`es que l'on fait un choix d'un sommet $s\in\BC\TC_0$ et une base d'un r\'eseau qui repr\'esente $s,$ l'isomorphisme dans le fait pr\'ec\'edent nous fournit un sous-groupe parabolique de $G_s/G_s^+,$ et donc un simplexe $\s$ contenant $s.$ De plus, on a un diagramme commutatif:
$$
  \begin{array}{ccc}
    G_s/G_s^+ & \cong & \GL_d(\FM_q) \\
    \bigcup &  & \bigcup \\
    G^+_\s/G^+_s & \cong & U_I \\
  \end{array}
$$

\begin{fait}
Sous l'hypoth\`ese du fait pr\'ec\'edent, l'isomorphisme $\o{\O}_s\simto \o{\O}^{d-1}_{\oFq}$ comme dans {\em loc. cit.} identifie $\o{\O}_\s$ \`a $\o{C}_I,$ o\`u $\o C_I:=(\o\O^{d-1}_{\FM_q})^{U_I}$ les points stables sous $U_I.$

\end{fait}
\begin{preuve}
Supposons tout d'abord que $I=\D\ba \{\a_i\}$ pour une racine simple $\a_i\in\D.$ D'apr\`es ce qui pr\'ec\`ede, $I$ fait un choix d'un simplexe $[s,s']$ contenant $s,$ o\`u $s'$ est un sommet voisin de $s.$ La remarque dans \cite[(4.1.4)]{wang_DL} nous dit que $\o{C}_I$ s'identifie \`a $\o{\O}_s\cap\o{\O}_{s'}=\o{\O}_{[s,s']}.$ Pour $I$ quelconque, supposons que $I=\D\ba\{\a_{i_1},\ldots,\a_{i_r}\}.$ Alors $I=I_1\cap\cdots\cap I_r,$ o\`u $I_j=\D\ba\{\a_{i_j}\},$ et $\o{C}_I=\o{C}_{I_1}\cap\cdots\cap \o{C}_{I_r}.$ Chaque $I_j$ fait un choix d'un sommet voisin $s_j$ de $s,$ et identifie $\o{C}_{I_j}$ \`a $\o{\O}_{[s,s_j]}.$ Notons que $\s:=[s,s_1,\ldots,s_r]$ est le simplexe donn\'e par $I.$ Par d\'efinition, $\o{\O}_\s=\o{\O}_{[s,s_1]}\cap\cdots\cap\o{\O}_{[s,s_r]},$ donc il s'identifie \`a $\o{C}_I.$
\end{preuve}

\end{rem}

\ali Rappelons le th\'eor\`eme principal (Th\'eor\`eme 5.1.1) de \cite{wang_DL}. Consid\'erons le diagramme suivant:
$$\xymatrix{
\DL^{d-1}_{\oFq}\ar[d]^{\pi} & \\
\O^{d-1}_{\oFq}\ar@{^(->}[r]^j & \o{\O}^{d-1}_{\oFq} \ar@{<-^)}[r]^{i_I} & \o{C}_I
}
$$

\begin{thm}
Le morphisme de restriction
$$
R\G(\O^{d-1}_{\oFq},\pi_*\L)=R\G(\o{\O}^{d-1}_{\oFq} ,Rj_*(\pi_*\L))\To{\res.}R\G(\o{C}_I,i^*_IRj_*(\pi_*\L))
$$
induit un isomorphisme
$$
R\G(\O^{d-1}_{\oFq},\pi_*\L)^{U_I}\simto R\G(\o{C}_I,i^*_IRj_*(\pi_*\L)).
$$
\end{thm}

Cela nous donne la cons\'equence suivante sur la cohomologie de $\Sig^{ca}.$
\begin{theo}\label{new7}
Soient $s$ un sommet de $\BC\TC$ et $\s$ un simplexe contenant $s$, alors le morphisme canonique $H^q_c(\nu^{-1}(|\s|^*),\Lambda)\to H^q_c(\nu^{-1}(|s|^*),\Lambda)$ provenant de l'immersion ouverte $\nu^{-1}(|\s|^*)\into\nu^{-1}(|s|^*)$ induit un isomorphisme
\ini\begin{equation}\label{equa1}
H^q_c(\nu^{-1}(|\s|^*),\Lambda)\simto H^q_c(\nu^{-1}(|s|^*),\Lambda)^{G_\s^+}.
\end{equation}
\end{theo}
\begin{preuve}
Par la dualit\'e de Poincar\'e, il suffit de montrer que le morphisme de restriction $R\G(\nu^{-1}(|s|^*),\L)\To{} R\G(\nu^{-1}(|\s|^*),\L)$ induit un ismorphisme $R\G(\nu^{-1}(|s|^*),\L)_{G_\s^+}\simto R\G(\nu^{-1}(|\s|^*),\L).$ D'apr\`es Berkovich \cite[Corollary 3.5]{berk-vanishing}, $$R\G(\nu^{-1}(|\s|^*),\L)=R\G(\o{\O}_\s,R\Psi_\eta(p_*\L)|_{\o{\O}_\s}),$$ o\`u $p$ est la projection $\Sig^{ca}\onto \O^{d-1,ca}_K.$

Notons $\o{p}_s:\o{\Sig}^0_s\to \o{\O}^0_s$ et $i_\s:\o{\O}_\s\into\o{\O}_s.$ D'apr\`es le th\'eor\`eme \ref{new2} et le corollaire \ref{new4}, on a un diagramme commutatif:
$$
\xymatrix{
R\G(\nu^{-1}(|s|^*),\L) \ar[d]^{\res.}\ar[r]^\cong & R\G(\o{\O}^0_s,\o{p}_{s*} \L)\ar[d]^{\res.} \\
R\G(\nu^{-1}(|\s|^*),\L) \ar[r]^\cong & R\G(\o{\O}_\s,i^*_\s Rj_{s*}\o{p}_{s*} \L)
}
$$
Donc on se ram\`ene \`a montrer que le morphisme de restriction
induit un isomorphisme
\begin{equation*}
R\G(\o{\O}^0_s,\o{p}_{s*}\L)_{G^+_\s/G_s^+}\simto R\G(\o{\O}_\s,i_\s^*Rj_{s,*}\o{p}_{s*}\L).
\end{equation*}
Comme $G^+_\s/G_s^+$ est un $p$-groupe fini et $\L$ est de torsion premier \`a $p$, on peut identifier canoniquement $R\G(\o{\O}^0_s,\o{p}_{s*}\L)_{G^+_\s/G_s^+}$ \`a $R\G(\o{\O}^0_s,\o{p}_{s*}\L)^{G^+_\s/G_s^+}.$ On obtient l'\'enonc\'e du th\'eor\`eme, en vertu du th\'eor\`eme pr\'ec\'edent et les isomorphismes canoniques:
$$
\xymatrix{
R\G(\o{\O}^0_s,\o{p}_{s*} \L)\ar[d]^{\res.} \ar[r]^\cong & R\G(\O^{d-1}_{\oFq},\pi_*\L)\ar[d]^{\res.}\\
 R\G(\o{\O}_\s,i^*_\s Rj_{s*}\o{p}_{s*} \L) \ar[r]^\cong & R\G(\o{C}_I,i^*_I Rj_{I*}\pi_*\L).
}
$$
\end{preuve}

\section{La partie supercuspidale de la cohomologie}

Dans cette section, on donne quelques cons\'equences sur la partie supercuspidale de la cohomologie de $\Sig^{ca}.$

\subsection{La d\'emonstration de Th\'eor\`eme A.}
Fixons $\ell\neq p$ un nombre premier et une cl\^oture alg\'ebrique $\oQl$ de $\QM_\ell.$

\ali\label{b} Soient $\s'\subset \s$ deux simplexes de $\BC\TC$, l'immersion ouverte $\nu^{-1}(|\s|^*)\into\nu^{-1}(|\s'|^*)$ induit un morphisme canonique:
\[
H^q_c(\nu^{-1}(|\s|^*),\Lambda)\longrightarrow H^q_c(\nu^{-1}(|\s'|^*),\Lambda).
\]
Donc les donn\'ees
$$
\begin{cases}
(\s\in\BC\TC)\mapsto H^q_c(\nu^{-1}(|\s|^*),\Lambda)\\
(\s'\subset \s)\mapsto \big(H^q_c(\nu^{-1}(|\s|^*),\Lambda)\to H^q_c(\nu^{-1}(|\s'|^*),\Lambda)\big)
\end{cases}
$$
d\'efinissent un syst\`eme de coefficients $G^\circ$-\'equivariant (voir \cite[3.2.3]{dat-drinfeld}) \`a valeurs dans les $\Lambda$-modules que nous noterons simplement $\s\mapsto H^q_c(\nu^{-1}(|\s|^*),\Lambda).$ Nous avons le fait suivant bien connu (voir \cite[Prop. 3.2.4]{dat-drinfeld} pour le fait et \cite[3.2.3]{dat-drinfeld} pour les notations):

\begin{fait}\label{suite spec}
Il existe une suite spectrale $G^\circ$-\'equivariante
\ini\begin{equation}\label{suite spec1}
E^{pq}_1=\CC^{or}_c(\BC\TC_{(-p)},\s\mapsto H^q_c(\nu^{-1}(|\s|^*),\Lambda))\Longrightarrow H^{p+q}_c(\Sigma^{ca},\Lambda)
\end{equation}
dont la diff\'erentielle $d^{pq}_1$ est celle du complexe de cha\^ines du syst\`eme de coefficients $\s\mapsto H^q_c(\nu^{-1}(|\s|^*),\Lambda).$
\end{fait}

\ali Notons $W_K$ le groupe de Weil associ\'e \`a $K$ et $I_K$ le sous-groupe d'inertie. Posons $GDW:=G\times D^\times \times W_K,$ et $v:GDW\to \ZM$ l'homorphisme qui envoie un \'el\'ement $(g,\d,w)$ vers l'entier $\val_K(\det(g^{-1})\Nr(\d)\Art^{-1}(w))\in\ZM,$ o\`u $\Nr:D^\times\to K^\times$ d\'esigne la norme r\'eduite et $\Art^{-1}:W_K\onto W_K^{ab}\to K^\times$ d\'esigne la compos\'ee de l'inverse du morphisme d'Artin qui envoie l'uniformisante $\varpi$ vers un Frobenius g\'eom\'etrique fix\'e $\varphi.$ On d\'esigne $(GDW)^0:=v^{-1}(0)$ le noyau de $v,$ et $[GDW]_{d}$ le sous-groupe distingu\'e form\'e des \'el\'ements $(g,\d,w)$ tels que $v(g,\d,w)\in d\ZM.$ On consid\`ere $G$ (resp. $D^\times,$ $W_K$) comme un sous-groupe de $GDW$ via l'inclusion naturelle, et on notera $[G]_{d}$ (resp. $[D]_{d},$ $[W_K]_{d}$) son intersection avec $[GDW]_{d}.$ D\`es que l'on identifie $K^\times$ au centre de $G,$ on a $[GDW]_d=(GDW)^0\varpi^\ZM.$

\ali On a vu dans le paragraphe \ref{o} que l'espace analytique $\Sig^{ca}$ admet une action naturelle du groupe $G^\circ\times\OC_D^\times\times I_K.$ Afin de d\'efinir une action du groupe $GDW,$ on suit la construction de Rapoport-Zink \cite{RZ} en adoptant les notations de \cite[\S 3]{Dat-elliptic}, i.e. on ne fixe pas la hauteur des quasi-isog\'enies qui rigidifient les probl\`emes de modules dans la d\'efinition du foncteur $G^{Dr}$ dans \ref{new1}. Ce nouveau probl\`eme de modules $\wt G$ est de m\^eme pro-repr\'esentable par un sch\'ema formel $\wh \MC_{Dr,0}$ sur $\Spf \breve{\OC},$ {\em cf.} \cite[(3.1.3)]{Dat-elliptic}. Gr\^ace \`a la d\'ecomposition suivant la hauteur de quasi-isog\'enie, $\wt G=\coprod_{h\in\ZM} G^{(h)}$ o\`u $G^{(h)}$ classifie les classes de triples $(\psi,X,\rho)$ avec $\rho$ de hauteur $dh.$ Chaque $G^{(h)}$ est isomorphe (non-canoniquement) \`a $G^{(0)}$ et $G^{(0)}$ est le foncteur $G^{Dr}$ de Drinfeld. On en d\'eduit un isomorphisme non-canonique $\wh \MC_{Dr,0}\cong \wh \MC_{Dr,0}^{(0)}\times\ZM$ o\`u $\wh \MC_{Dr,0}^{(0)}:=\wh \O^{d-1}_{\breve\OC}$ est le sch\'ema formel qui repr\'esente $G^{Dr}.$ On note $\wt \XG$ le $\OC_D$-module formel sp\'ecial universel sur $\wh \MC_{Dr,0}.$ Le noyau $\wt\XG[\Pi_D]$ de la multiplication par $\Pi_D$ dans $\wt\XG$ est un sch\'ema formel en groupe fini plat de rang $q^d$ au-dessus de $\wh \MC_{Dr,0},$ et qui est \'etale en fibre g\'en\'erique. Le $(\OC_D/\Pi_D\OC_D)^\times$-torseur $\underline{\Isom}_{\OC_D}(\OC_D/\Pi_D\OC_D,\wt\XG[\Pi_D]^{an})$ sur $\MC_{Dr,0}$ la fibre g\'en\'erique au sens de Raynaud-Berkovich de $\wh \MC_{Dr,0},$ est donc repr\'esent\'e par un $\breve{K}$-espace analytique $\MC_{Dr,1},$ qui est un rev\^etement \'etale de $\MC_{Dr,0}$ de groupe $\OC_D^\times/(1+\Pi_D\OC_D)\cong \FM_{q^d}^\times.$ En posant $\MC_{Dr,1}^{ca}:=\MC_{Dr,1}\wh\otimes_{\breve{K}}\wh{K^{ca}},$ on a un isomorphisme $\MC_{Dr,1}^{ca}\cong\Sig^{ca}\times\ZM.$ Suivant \cite[(3.1.5)]{Dat-elliptic}, on d\'efinit une action du groupe $GDW$ sur $\MC_{Dr,1}^{ca}$ telle que la composante $\Sig^{ca}\times\{0\}$ soit stable sous l'action du sous-groupe $(GDW)^0.$ D\'esormais, on identifie $\Sig^{ca}$ au sous-espace analytique ouvert-ferm\'e $\Sig^{ca}\times\{0\}$ de $\MC_{Dr,1}^{ca}.$

\ali\label{Sec:6} Soit $\th:\FM_{q^d}^\times\to\oQl^\times$ un caract\`ere $d$-primitif, i.e. $\th$ ne se factorise pas par la norme $N_{\FM_{q^d}/\FM_{q^{f}}}:\FM_{q^{d}}^\times\onto\FM_{q^{f}}^\times,~\forall f|d$ et $f\neq d.$ D\'efinissons les repr\'esentations suivantes:

\begin{itemize}
  \item une repr\'esentation $\rho(\th):=\ind_{[D]_d}^{D^\times}\th$ de $D^\times,$ o\`u $[D]_d=\OC_D^\times\varpi^\ZM$ tel que $\OC^\times_D$ agit via la projection $\OC_D^\times\onto \FM^\times_{q^d},$ et $\varpi^\ZM$ agit trivialement.
  \item une repr\'esentation $\pi(\th):=\ind_{\GL_d(\OC)\varpi^\ZM}^{G}\o\pi_\th$ de $G,$ o\`u $\o\pi_\th$ est la repr\'esentation cuspidale de $\GL_d(\FM_q)$ associ\'ee \`a $\th$ via la correspondance de Green (ou de Deligne-Lusztig), vue comme une repr\'esentation de $\GL_d(\OC)\varpi^\ZM$ telle que $\GL_d(\OC)$ agit via la projection $\GL_d(\OC)\onto\GL_d(\FM_q)$ et $\varpi^\ZM$ agit trivialement.
  \item une repr\'esentation $\s^\sharp(\th):=\ind_{[W_K]_d}^{W_K}\wt\th$ de $W_K,$ o\`u $\wt\th$ est le caract\`ere de $[W_K]_d=I_K\langle\varphi^d\rangle^\ZM$ tel que $\wt\th(\varphi^d)=(-1)^{d-1}q^{\frac{d(d-1)}{2}},$ et $\wt\th|_{I_K}$ se factorise par l'inertie mod\'er\'ee $I_K\onto I_K/I_{K(\varpi_t)}\cong\mu_{q^d-1}\To{\th}\oQl^\times,$ ici $\varpi_t$ est la racine $(q^d-1)$-i\`eme de $\varpi$ fix\'ee dans \ref{Sec:4}.
\end{itemize}

\medskip
Posons
$$
\Hb^{i}_c:=H^i_c(\MC_{Dr,1}^{ca}/\varpi^\ZM,\oQl), ~i\in\{d-1,\ldots, 2d-2\}.
$$
D'apr\`es la description pr\'ec\'edente,
$$
\Hb^i_c=\ind^{GDW}_{[GDW]_d}H^{i}_c(\Sig^{ca},\oQl),
$$
o\`u $\varpi\in K^\times\subset G$ agit trivialement. On consid\`ere la partie $\rho(\th)$-isotypique de $\Hb^{i}_c.$

\begin{theo}\label{T:1}
Pour tout $\th:\FM_{q^d}^\times\to\oQl^\times$ caract\`ere $d$-primitif, on a

$$
\Hom_{D^\times}(\rho(\th),\Hb^{i}_c)\underset{G\times W_K}{\simeq}\left\{
                                                          \begin{array}{ll}
                                                            \pi(\th)\otimes\s^\sharp(\th), & \hbox{si $i=d-1$;} \\
                                                            0, & \hbox{si $i\neq d-1$.}
                                                          \end{array}
                                                        \right.
$$

\end{theo}
\begin{preuve}
Ceci d\'ecoule du lemme \ref{L:1} et la proposition \ref{Prop3}.
\end{preuve}

\begin{lemme}\label{L:2}En tant que repr\'esentation de $GW:=G\times W_K,$ on a
$$
\Hom_{D^\times}(\rho(\th),\Hb^{i}_c)=\ind^{GW}_{[GW]_d}\Hom_{\FM^\times_{q^d}}(\th,H^{i}_c(\Sig^{ca},\oQl)),
$$
o\`u $[GW]_d=\{(g,w)\in GW~|~v(g,1,w)\in d\ZM\}.$
\end{lemme}
\begin{preuve}
D'apr\`es la r\'eciprocit\'e de Frobenius, on a
\begin{align*}
\Hom_{D^\times}(\rho(\th),\Hb^{i}_c)&=\Hom_{[D]_d}(\th,\ind^{GDW}_{[GDW]_d}H^{i}_c(\Sig^{ca},\oQl))\\
&=\ind^{GW}_{[GW]_d}\Hom_{\FM^\times_{q^d}}(\th,H^{i}_c(\Sig^{ca},\oQl)).
\end{align*}
\end{preuve}

\begin{lemme}\label{L:1}
En tant que repr\'esentation de $G,$ on a
$$
\Hom_{D^\times}(\rho(\th),\Hb^{i}_c)=\left\{
                                                          \begin{array}{ll}
                                                            \pi(\th)^{\oplus d}, & \hbox{si $i=d-1$;} \\
                                                            0, & \hbox{sinon.}
                                                          \end{array}
                                                        \right.
$$
\end{lemme}
\begin{preuve}
Notons $M(\th)$ la partie $\th$-isotypique d'un $\FM_{q^d}^\times$-module $M.$ Rappelons que si $s$ est un sommet de $\BC\TC,$ nous avons pour tout $i\in\NM$ un isomorphisme ({\em cf.} \ref{C1C1})
\[
H^i_c(\nu^{-1}(|s|^*),\oQl)\cong H^i_c(\DL^{d-1}_{\oFq},\oQl).
\]
Comme $\th$ est $d$-primitif, gr\^ace \`a la propri\'et\'e des vari\'et\'es de Deligne-Lusztig \cite[Lemma 9.14]{deligne-lusztig}, on a
$$
H^i_c(\DL^{d-1}_{\oFq},\oQl)(\th)\cong\left\{
                                       \begin{array}{ll}
                                         \o\pi_{\th}, & \hbox{si $i=d-1$;} \\
                                         0, & \hbox{sinon.}
                                       \end{array}
                                     \right.
$$
En vertu du th\'eor\`eme \ref{new7}, si $i\neq d-1,$ $H^i_c(\nu^{-1}(|\s|^*),\oQl)(\th)=0 ~\forall \s.$ Si $i=d-1,$ $\s\ni s$ un simplexe tel que $\dim\s>0,$ comme $\o\pi_{\th}$ est une repr\'esentation cuspidale de $\GL_d(\FM_q),$ on a
$$
H^{d-1}_c(\nu^{-1}(|\s|^*),\oQl)(\th)\cong \o\pi_{\th}^{G^+_\s/G^+_s}=0.
$$
Donc, la partie $\th$-isotypique de la suite spectrale \ref{suite spec1} d\'eg\'en\`ere en $E_1$ lorsque $i=d-1.$ En conclusion, on a
\ini\begin{align}\label{R1E1}
\Hom_{\FM^\times_{q^d}}(\th,H^{i}_c(\Sig^{ca},\oQl))=\left\{
                                       \begin{array}{ll}
                                         \bigoplus_{s\in\BC\TC_0}H^{d-1}_c(\nu^{-1}(|s|^*),\oQl)(\theta)=\bigoplus_{s\in\BC\TC_0}\o\pi_\th, & \hbox{si $i=d-1$;} \\
                                         0, & \hbox{sinon.}
                                       \end{array}
                                     \right.
\end{align}
Ceci nous fournit un isomorphisme de $\oQl[G]_d$-modules
\begin{equation*}
\Hom_{\FM^\times_{q^d}}(\th,H^{d-1}_c(\Sig^{ca},\oQl))\cong \pi(\th)|_{[G]_d}.
\end{equation*}
Alors, en tant que repr\'esentation de $G,$ nous avons
\begin{align*}
\Hom_{D^\times}(\rho(\th),\Hb^{d-1}_c)&=\ind^{G}_{[G]_d}\Hom_{\FM^\times_{q^d}}(\th,H^{d-1}_c(\Sig^{ca},\oQl))\\
&=\ind^G_{[G]_d}\pi(\th)|_{[G]_d}
\end{align*}
Notons $\wt\pi(\th):=\ind^{[G]_d}_{\GL_d(\OC)\varpi^\ZM}\o\pi_\th,$ il s'ensuit que
$$
\ind^G_{[G]_d}\pi(\th)|_{[G]_d}=\ind^G_{[G]_d}(\bigoplus_{x\in G/[G]_d}~ ^x\wt\pi(\th))=\bigoplus_{x\in G/[G]_d}\ind^G_{[G]_d} (^x\wt\pi(\th)),
$$
o\`u $^x\wt\pi(\th)(g)=\wt\pi(\th)(x^{-1}gx)$ pour $g\in[G]_d.$ D'apr\`es la formule de Mackey, on sait que
$$
\ind^G_{[G]_d} (^x\wt\pi(\th))=\pi(\th),~\forall x\in G/[G]_d.
$$
Ceci nous permet de conclure car $|G/[G]_d|=d.$
\end{preuve}

Consid\'erons le $W_K$-module $$\Hom_G(\pi(\th),\Hom_{D^\times}(\rho(\th),\Hb^{d-1}_c)).$$
Dans le reste de ce paragraphe, on d\'emontre la proposition suivante:

\begin{prop}\label{Prop3}
On a un isomorphisme de $W_K$-modules
$$
\Hom_G(\pi(\th),\Hom_{D^\times}(\rho(\th),\Hb^{d-1}_c))\cong \s^\sharp(\th).
$$
\end{prop}
\begin{preuve}
D'apr\`es la r\'eciprocit\'e de Frobenius, on a
\begin{align*}
\Hom_G(\pi(\th),\Hom_{D^\times}(\rho(\th),\Hb^{d-1}_c))&=\Hom_{G}(\ind^G_{\GL_d(\OC)\varpi^\ZM}\o\pi_{\th},\ind^{GW}_{[GW]_d}\Hom_{\FM^\times_{q^d}}(\th,H^{d-1}_c(\Sig^{ca},\oQl))\\
&=\ind^{W_K}_{[W_K]_d}\Hom_{\GL_d(\OC)\varpi^{\ZM}}(\o\pi_\th,\Hom_{\FM^\times_{q^d}}(\th,H^{d-1}_c(\Sig^{ca},\oQl))\\
\end{align*}
D'apr\`es \ref{R1E1}, 
\begin{align*}
\Hom_G(\pi(\th),\Hom_{D^\times}(\rho(\th),\Hb^{d-1}_c))&=\ind^{W_K}_{[W_K]_d}\Hom_{\GL_d(\OC)\varpi^{\ZM}}\big(\o\pi_\th,\bigoplus_{s'\in\BC\TC_0}H^{d-1}_c(\nu^{-1}(|s'|^*),\oQl)(\th)\big)\\
&=\ind^{W_K}_{[W_K]_d}\Hom_{\GL_d(\FM_q)}\Big(\o\pi_\th,\big(\bigoplus_{s'\in\BC\TC_0}H^{d-1}_c(\nu^{-1}(|s'|^*),\oQl)(\th)\big)^{G^+_s}\Big),
\end{align*}
o\`u $s$ d\'esigne le r\'eseau standard $[\OC^d].$
On d\'emontre que 
$$
\big(\bigoplus_{s'\in\BC\TC_0}H^{d-1}_c(\nu^{-1}(|s'|^*),\oQl)(\th)\big)^{G^+_s}=H^{d-1}_c(\nu^{-1}(|s|^*),\oQl)(\th).
$$
\'Evidemment, les \'el\'ements dans $H^{d-1}_c(\nu^{-1}(|s|^*),\oQl)(\th)$ sont $G^+_s$-invariants. Soit $x=\sum x_{s'}\in\big(\bigoplus_{s'\in\BC\TC_0}H^{d-1}_c(\nu^{-1}(|s'|^*),\oQl)(\th)\big)^{G^+_s}$ avec $x_{s'}\in H^{d-1}_c(\nu^{-1}(|s'|^*),\oQl)(\th).$ Notons $e_{G^+_{s'}}$ l'idempotent central associ\'e \`a pro-$p$-groupe $G^+_{s'}.$ Par hypoth\`ese, on a $e_{G^+_s}x=x$ et $e_{G^+_{s'}}x_{s'}=x_{s'},~\forall s'.$ Soit $s'$ un sommet diff\'erent de $s,$ on a 
$e_{G^+_{s}}x_{s'}=0.$ En effet, il existe un sommet $s''$ dans l'enclos de $s$ et $s'$ tel que $s', s''$ soient adjacent. D'apr\`es \cite[2.2]{MS}, $e_{G^+_{s}}e_{G^+_{s'}}=e_{G^+_{s}}e_{G^+_{s''}}e_{G^+_{s'}}.$ Notons que $G^+_{s''}G^+_{s'}/G^+_{s'}$ est isomorphe \`a un sous-groupe unipotent non trivial de $\GL_d(\FM_q),$ et que $H^{d-1}_c(\nu^{-1}(|s'|^*),\oQl)(\th)$ est cuspidale en tant que repr\'esentation de $G_{s'}/G^+_{s'}\cong\GL_d(\FM_q).$ On a alors
$$
e_{G^+_{s}}x_{s'}=e_{G^+_{s}}e_{G^+_{s'}}x_{s'}=e_{G^+_{s}}(e_{G^+_{s''}}e_{G^+_{s'}}x)=0.
$$ 
On en d\'eduit que
$$
x=e_{G^+_s}x=e_{G^+_s}x_s+\sum_{s'\neq s}e_{G^+_s}x_{s'}=x_s\in H^{d-1}_c(\nu^{-1}(|s|^*),\oQl)(\th).
$$
Par cons\'equent, on a
$$
\Hom_G(\pi(\th),\Hom_{D^\times}(\rho(\th),\Hb^{d-1}_c))=\ind^{W_K}_{[W_K]_d}\Hom_{\GL_d(\FM_q)}(\o\pi_\th,H^{d-1}_c(\nu^{-1}(|s|^*),\oQl)(\th)).
$$
Notons que $$V:=\Hom_{\GL_d(\FM_q)}(\o\pi_{\th},H^{d-1}_c(\nu^{-1}(|s|^*),\oQl)(\th))$$ est un $\oQl$-espace vectoriel de dimension $1,$ car $$H^{d-1}_c(\nu^{-1}(|s|^*),\oQl)(\th)\cong H^{d-1}_c(\DL_{\oFq},\oQl)(\th)=\o\pi_{\th}$$ en tant que repr\'esentation de $\GL_d(\FM_q).$ Donc on se ram\`ene \`a comprendre l'action de $[W_K]_d=I_K\langle\varphi^d\rangle^\ZM$ sur $V.$

Tout d'abord, on note que la dualit\'e de Poincar\'e entra\^ine que 
\begin{align*}
V&=\Hom_{\GL_d(\FM_q)}(\o\pi_{\th},H^{d-1}(\nu^{-1}(|s|),\oQl(d-1))(\th^{-1})^\vee).\\
&=\Hom_{\GL_d(\FM_q)}(\o\pi_{\th},H^{d-1}(\Sig^{ca}_s,\oQl(d-1))(\th^{-1})^\vee).
\end{align*}
En comparant \ref{Sec:4} et \ref{Sec:1}, on observe que l'action de $I_K$ sur $\wh\Sig^0_s$ ({\em cf.} \ref{Sec:4}) se factorise par $I_K/I_{K(\varpi_t)}\cong\FM^\times_{q^d},$ et elle s'identifie \`a l'action de $\FM^\times_{q^d}$ sur $\DL^{d-1}_{\oFq}.$ Alors l'action de $I_K$ sur $V$ se factorise par $I_K/I_{K(\varpi_t)}$ via le caract\`ere $\th.$

On est amen\'e donc \`a \'etudier l'action de $\varphi^d$ sur $V.$ \`A priori, l'espace analytique $\Sig^{ca}$ n'a qu'une action du groupe d'inertie $I_K.$  Cependant la donn\'ee de descente \`a la Weil sur $\wh\MC_{Dr,0}$ (\cite[(3.48)]{RZ}) d\'efinit une action de l'\'el\'ement de Frobenius g\'eom\'etrique $\varphi$ sur $\MC_{Dr,1}^{ca}.$ Pour d\'efinir une donn\'ee de descente de $\Sig,$ il faut d\'ecaler par l'endomorphisme $\Pi_D^{-1}.$ Pr\'ecis\'ement, d'apr\`es \cite[(3.72)]{RZ}, le morphisme $\Pi_D^{-1}\circ\varphi$ induit la donn\'ee de descente canonique de $\wh\MC^{(0)}_{Dr,0},$ donc il induit une donn\'ee de descente de $\Sig,$ et de m\^eme une donn\'ee de descente de $\wh\Sig_s^0.$ En particulier, $\Pi_D^{-1}\circ\varphi$ induit une donn\'ee de descente de la fibre sp\'eciale $\o\Sig_s^0$ de $\wh\Sig_s^0.$ Autrement dit, cela nous fournit une structure $\FM_q$-rationnelle de $\o\Sig_s^0,$ et on note $\Fr$ le morphisme de Frobenius associ\'e.

Pour d\'eterminer cette structure $\FM_q$-rationnelle, on consid\`ere l'action de $\Fr$ sur les composantes connexes g\'eom\'etriques de $\o\Sig_s^0.$ Notons $\pi_0(X)$ l'ensemble des composantes connexes g\'eom\'etriques de $X$ pour un espace analytique ou une vari\'et\'e $X.$  On a le fait suivant:

\begin{fait}
$\pi_0(\MC_{Dr,1})$ est un espace homog\`ene sous $K^\times,$ isomorphe \`a $K^\times/(1+\varpi\OC),$ sur lequel
l'action de $G$ est donn\'ee par $G\To{\det}K^\times,$
celle de $D^\times$ est donn\'ee par $D^\times\To{\Nr}K^\times,$
celle de $W_K$ est donn\'ee par $\Art^{-1}_K:W_K\onto W_K^{ab}\simto K^\times.$

\end{fait}

\begin{rem}
Dans \cite{chen}, Chen a \'etudi\'e les composantes connexes g\'eom\'etriques de la tour de l'espace de Rapoport-Zink non-ramifi\'e de certains types, avec les actions de diff\'erents groupes associ\'es \`a la donn\'ee de Rapoport-Zink. En particulier, elle a d\'ecrit les composantes connexes g\'eom\'etriques de la tour de Lubin-Tate $(\MC_{LT,n})_{n\in\NM}$ munies des actions de $G,D^\times,W_K.$ En vertu de l'isomorphisme de Faltings-Fargues \cite{FGL} entre la tour de Lubin-Tate et la tour de Drinfeld et en prenant le sous-ensemble fix\'e par $1+\Pi_D\OC_D,$ on obtient l'\'enonc\'e du fait pr\'ec\'edent.
\end{rem}

\begin{lem}
$\pi_0(\o\Sig^0_s)$ est un espace principal homog\`ene sous $\FM_q^\times$ sur lequel $\Fr$ agit par multiplication par $(-1)^{d-1}.$
\end{lem}

\begin{preuve}
Notons que $\Nr(\Pi_D)=(-1)^{d-1}\varpi$ et $\Art^{-1}(\ph)=\varpi.$ D'apr\`es le fait pr\'ec\'edent, $\pi_0(\Sig)$ est un espace principal homog\`ene sous $\FM_q^\times,$ sur lequel $\Pi_D^{-1}\circ\varphi$ agit par multiplication par  $$\Nr(\Pi_D^{-1})\cdot\Art^{-1}(\varphi)=(-1)^{d-1}\varpi^{-1}\cdot\varpi=(-1)^{d-1}.$$
Comme $\wh\Sig_s^0$ est un $\mu_{q^d-1}$-torseur sur $\wh\O^{d-1}_{\breve{\OC},s},$ on a une bijection canonique entre $\pi_0(\nu^{-1}(|s|))$ et $\pi_0(\o\Sig_s^0).$ D'apr\`es l'\'equation \ref{Eq:2} de $\o\Sig_s^0\cong \DL^{d-1}_{\oFq},$ on sait que $|\pi_0(\o\Sig_s^0)|=q-1.$ Les immersions $\nu^{-1}(|s|)\to \nu^{-1}(|s|^*)\to\Sig$ induisent les morphismes $\a,\b$ entre les composantes connexes
$$
\pi_0(\nu^{-1}(|s|))\To{\a}\pi_0(\nu^{-1}(|s|^*))\To{\b}\pi_0(\Sig).
$$
On d\'emontre que le compos\'e $\b\circ\a$ est une bijection. En effet, $\a$ est bijective, car $H^0(\nu^{-1}(|s|),\L)=H^0(\nu^{-1}(|s|^*),\L),$ d'apr\`es \ref{new2}. Puisque $|\pi_0(\nu^{-1}(|s|))|=|\pi_0(\Sig)|=q-1,$ il suffit de montrer que $\b$ est surjective. Notons que $\{\nu^{-1}(|s'|^*)\}_{s'\in\BC\TC_0}$ forment un recourement ouvert de $\Sig.$ On se ram\`ene alors \`a prouver que pour $s'$ un sommet quelconque et $U'$ une composante connexe de $\pi_0(\nu^{-1}(|s|^*))$, il existe une suite de sommets $s_0=s, s_1,\ldots,s_n=s'$ et des composantes conexes $U_i\in\pi_0(\nu^{-1}(|s_i|^*))$ avec $U_n=U'$ telles que $s_i,s_{i+1}$ soient adjacent et $U_i\cap U_{i+1}\neq\emptyset,~\forall i.$ Par d\'evissage, il suffit de consid\'erer le cas o\`u $s$ et $s'$ sont adjacent. Notons $\s$ le simplexe $[s,s'],$ et consid\'erons les immersions ouvertes $\nu^{-1}(|\s|^*)\into\nu^{-1}(|s|^*)$ et $\nu^{-1}(|\s|^*)\into\nu^{-1}(|s'|^*)$ induisant les morphismes entre $\pi_0:$
$$
\a_1:\pi_0(\nu^{-1}(|\s|^*))\to\pi_0(\nu^{-1}(|s|^*)),~\a_2:\pi_0(\nu^{-1}(|\s|^*))\to\pi_0(\nu^{-1}(|s'|^*)).
$$  
Notons que $\a_1$ et $\a_2$ sont bijectives. Ceci d\'ecoule du th\'eor\`eme \ref{new7} et du fait que le $p$-groupe $G^+_\s$ agit trivialement sur $H^{2d-2}_c(\nu^{-1}(|s|^*),\L).$ Donc, pour $U'\in\pi_0(\nu^{-1}(|s'|^*)),$ il existe une composante connexe $U\in\pi_0(\nu^{-1}(|s|^*))$ tel que $U\cap U'$ contienne une composante connexe du ouvert $\nu^{-1}(|\s|^*)$ de $\Sig.$ Par cons\'equent, $\b$ est bijective. 

On a alors une bijection canonique entre $\pi_0(\Sig)$ et $\pi_0(\o\Sig_s^0)$ compatible avec l'action de $\FM_q^\times.$ Par cons\'equent, le morphisme de Frobenius $\Fr$ induit la multiplication par $(-1)^{d-1}$ sur $\pi_0(\o\Sig_s^0). $
\end{preuve}

\begin{lemme}\label{4.2Lemme1}
La forme $\FM_q$-rationnelle de $\o\Sig_s^0$ induite par la donn\'ee de descente $\Fr$ est isomorphe \`a la vari\'et\'e de Deligne-Lusztig ($\FM_q$-rationnelle) $\DL^{d-1}$ d\'efinie par l'\'equation \ref{Eq:2}.
\end{lemme}
\begin{preuve}
D'apr\`es le th\'eor\`eme \ref{prop2}, on sait que $\o\Sig^0_s\cong\DL^{d-1}_{\oFq}.$ Gr\^ace \`a la suite spectrale
$$
E^{pq}_2=H^p(\Gal(\oFq/\FM_q),H^q(\O^{d-1}_{\oFq},\mu_{q^d-1}))\ddonc H^{p+q}(\O^{d-1}_{\FM_q},\mu_{q^d-1}),
$$
on a une suite exacte
$$
0\to H^1(\Gal(\oFq/\FM_q),\mu_{q^d-1})\to H^{1}(\O^{d-1}_{\FM_q},\mu_{q^d-1})\to H^1(\O^{d-1}_{\oFq},\mu_{q^d-1})\to H^2(\Gal(\oFq/\FM_q),\mu_{q^d-1}).
$$
Notons que $H^1(\Gal(\oFq/\FM_q),\mu_{q^d-1})=\FM_q^\times$ gr\^ace \`a la suite exacte de Kummer et Thm. 90 de Hilbert, et que $H^2(\Gal(\oFq/\FM_q),\mu_{q^d-1})=0.$ On en d\'eduit que les formes $\FM_q$-rationnelles de $\DL_{\oFq}$ sont donn\'ees par $Y_a$ pour tout $a\in\FM_q^\times,$ o\`u $Y_a$ est donn\'ee par l'\'equation $\det((X_i^{q^j})_{0\leq i,j\leq d-1})^{q-1}=a.$ En particulier, $\DL^{d-1}=Y_{(-1)^{d-1}}.$

L'ensemble des composantes connexes g\'eom\'etriques $\pi_0(Y_a)$ de $Y_a$ s'identifie \`a l'ensemble des racines $(q-1)$-i\`emes de $a.$ Pr\'ecis\'ement, on associe $\e$ telle que $\e^{q-1}=a$ la composante connexe g\'eom\'etrique d\'efinie par l'\'equation $\det((X_i^{q^j})_{0\leq i,j\leq d-1})=\e.$ Via cette identification, $t\in\FM_q^\times$ agit par multiplication qui envoie $\e$ vers $t\e.$ Consid\'erons le morphisme de Frobenius $\Frob:(X_i)_{0\leq i\leq d-1}\mapsto (X_i^q)_{0\leq i\leq d-1}$ sur $Y_a.$ Pour un point $x=(x_0,\ldots,x_{d-1})$ tel que $\det((x_i^{q^j})_{0\leq i,j\leq d-1})=\e,$ $\Frob(x)=(x_0^q,\ldots,x_{d-1}^q)$ satisfait l'\'equation $\det((x^q_i)^{q^j})_{0\leq i,j\leq d-1}=\e^q=a\cdot\e.$ Autrement dit, $\Frob$ agit sur $\pi_0(Y_a)$ par multiplication par $a.$ Il s'ensuit que $Y_{(-1)^{d-1}}=\DL^{d-1}$ est la forme $\FM_q$-rationnelle de $\o\Sig^0_s$ induite par $\Fr,$ en vertu du lemme pr\'ec\'edent.
\end{preuve}
\begin{theo}(Digne et Michel \cite{Digne-Michel})\label{4.2Thm1}
$\Fr^d$ agit sur $H_c^{d-1}(\o\Sig^0_s,\oQl)(\th)$ par multiplication par $(-1)^{d-1}q^{\frac{d(d-1)}{2}}.$
\end{theo}

\begin{preuve}
D'apr\`es le lemme pr\'ec\'edent, il suffit de montrer que $(-1)^{d-1}q^{\frac{d(d-1)}{2}}$ est l'unique valeur propre de $\Frob^d$ sur $H^{d-1}_c(\DL^{d-1},\oQl)(\th).$ Ceci est essentiellement d\'emontr\'e par Digne et Michel \cite[V. (3.14)]{Digne-Michel}. Cependant, ils ont calcul\'e la valeur propre de $\Frob^d$ pour la vari\'et\'e d\'efinie par l'\'equation $\det((X_i^{q^j})_{0\leq i,j\leq d-1})^{q-1}=1.$ De plus, dans leur formule (page 106 de {\em loc. cit.}), il manque un signe $(-1)^{d-1}$ ($d=n$ dans {\em loc. cit.}), comme ils ont oubli\'e qu'ils travaillaient avec une somme altern\'ee. En plus, le facteur $\th((-1)^{d-1})$ dans {\em loc. cit.} disparait. Cela vient du fait que l'on peut rendre un isomorphisme entre $Y_1$ et $Y_a$ par un changement de variables $f:X_i\mapsto X_i/b$ o\`u $b^{q^d-1}=a.$ Alors l'action de $\Frob^d$ sur $H^i_c(Y_a,\oQl)$ est la m\^eme que celle sur $H^i_c(Y_1,\oQl)$ tordue par l'automorphisme induit par le changement de variables $X_i\mapsto a\cdot X_i,$ donc n'est autre que l'action de l'\'el\'ement $a\in \FM_{q^d}^\times.$ Ce dernier agit bien par $\th(a)$ sur la partie $\th$-isotypique de la cohomologie. Par cons\'equent, dans notre cas o\`u $a=(-1)^{d-1},$ l'endomorphisme $\Frob^d$ agit sur $H^{d-1}_c(\DL^{d-1},\oQl)(\th)$ par la constante $(-1)^{d-1}q^{\frac{d(d-1)}{2}}.$
\end{preuve}
Revenons \`a la preuve de la proposition \ref{Prop3}. L'action de $\varphi^d$ sur $V$ est donn\'ee par celle sur $H^{d-1}(\Sig_s^{ca},\oQl(d-1))(\th)^\vee,$ donc elle s'identifie \`a l'action de $\Fr^d$ sur $H^{d-1}_c(\o\Sig^0_s,\oQl)(\th).$ En vertu du lemme pr\'ec\'edent, on obtient l'\'enonc\'e de la proposition.
\end{preuve}

\def\cprime{$'$} \def\cprime{$'$} \def\cprime{$'$} \def\cprime{$'$}
  \def\cprime{$'$} \def\cprime{$'$} \def\cprime{$'$}

\noindent
\textsc{Haoran Wang}\\
 Universit\'e Pierre et Marie Curie, Institut de Math\'ematiques de Jussieu

\medskip

\noindent{\it Adresse Pr\'esente:}
Max-Planck-Institut f\"ur Mathematik, Vivatsgasse 7, 53111 Bonn, Germany

\noindent\texttt{haoran.wang@mpim-bonn.mpg.de}


\begin{thebibliography}{Wan13b}
\expandafter\ifx\csname fonteauteurs\endcsname\relax
\def\fonteauteurs{\scshape}\fi

\bibitem[BC91]{boutot-carayol}
J.-F. \bgroup\fonteauteurs\bgroup Boutot\egroup\egroup{} et
  H.~\bgroup\fonteauteurs\bgroup Carayol\egroup\egroup{} :
\newblock Uniformisation {$p$}-adique des courbes de {S}himura: les
  th\'eor\`emes de \v {C}erednik et de {D}rinfel\cprime d.
\newblock {\em Ast\'erisque}, (196-197)\string:\penalty500\relax 7, 45--158
  (1992), 1991.
\newblock Courbes modulaires et courbes de Shimura (Orsay, 1987/1988).

\bibitem[Ber93]{berk-ihes}
Vladimir~G. \bgroup\fonteauteurs\bgroup Berkovich\egroup\egroup{} :
\newblock \'{E}tale cohomology for non-{A}rchimedean analytic spaces.
\newblock {\em Publ. Math. Inst. Hautes \'Etudes Sci.},
  78\string:\penalty500\relax 5--161, 1993.

\bibitem[Ber95]{berkovich-drin}
Vladimir~G. \bgroup\fonteauteurs\bgroup Berkovich\egroup\egroup{} :
\newblock The automorphism group of the {D}rinfel\cprime d half-plane.
\newblock {\em C. R. Acad. Sci. Paris S\'er. I Math.},
  321(9)\string:\penalty500\relax 1127--1132, 1995.

\bibitem[Ber96]{berk-vanishing}
Vladimir~G. \bgroup\fonteauteurs\bgroup Berkovich\egroup\egroup{} :
\newblock Vanishing cycles for formal schemes. {II}.
\newblock {\em Invent. Math.}, 125(2)\string:\penalty500\relax 367--390, 1996.

\bibitem[Boy99]{Boyer-these}
Pascal \bgroup\fonteauteurs\bgroup Boyer\egroup\egroup{} :
\newblock Mauvaise r\'eduction des vari\'et\'es de {D}rinfeld et correspondance
  de {L}anglands locale.
\newblock {\em Invent. Math.}, 138(3)\string:\penalty500\relax 573--629, 1999.

\bibitem[Car90]{Carayol-LT}
H.~\bgroup\fonteauteurs\bgroup Carayol\egroup\egroup{} :
\newblock Nonabelian {L}ubin-{T}ate theory.
\newblock \emph{In} {\em Automorphic forms, {S}himura varieties, and
  {$L$}-functions, {V}ol.\ {II} ({A}nn {A}rbor, {MI}, 1988)}, volume~11 de {\em
  Perspect. Math.}, pages 15--39. Academic Press, Boston, MA, 1990.

\bibitem[Che13]{chen}
Miaofen \bgroup\fonteauteurs\bgroup Chen\egroup\egroup{} :
\newblock {C}omposantes connexes g{\'e}om{\'e}triques de la tour des espaces de
  modules de groupes $p$-divisibles.
\newblock {\em \`A para\^itre dans Ann. Sci. \'Ecole Norm. Sup.}, 2013.

\bibitem[Dat06]{dat-drinfeld}
Jean-Fran{\c{c}}ois \bgroup\fonteauteurs\bgroup Dat\egroup\egroup{} :
\newblock Espaces sym\'etriques de {D}rinfeld et correspondance de {L}anglands
  locale.
\newblock {\em Ann. Sci. \'Ecole Norm. Sup. (4)},
  39(1)\string:\penalty500\relax 1--74, 2006.

\bibitem[Dat07]{Dat-elliptic}
Jean-Fran\c{c}ois \bgroup\fonteauteurs\bgroup Dat\egroup\egroup{} :
\newblock {T}h\'eorie de {L}ubin-{T}ate non-ab\'elienne et repr\'esentations
  elliptiques.
\newblock {\em Invent. Math.}, 169(1)\string:\penalty500\relax 75--152, 2007.

\bibitem[DH87]{DH}
Pierre \bgroup\fonteauteurs\bgroup Deligne\egroup\egroup{} et Dale
  \bgroup\fonteauteurs\bgroup Husemoller\egroup\egroup{} :
\newblock Survey of {D}rinfel\cprime d modules.
\newblock \emph{In} {\em Current trends in arithmetical algebraic geometry
  ({A}rcata, {C}alif., 1985)}, volume~67 de {\em Contemp. Math.}, pages 25--91.
  Amer. Math. Soc., Providence, RI, 1987.

\bibitem[DL76]{deligne-lusztig}
P.~\bgroup\fonteauteurs\bgroup Deligne\egroup\egroup{} et
  G.~\bgroup\fonteauteurs\bgroup Lusztig\egroup\egroup{} :
\newblock Representations of reductive groups over finite fields.
\newblock {\em Ann. of Math. (2)}, 103(1)\string:\penalty500\relax 103--161,
  1976.

\bibitem[DM85]{Digne-Michel}
Fran{\c{c}}ois \bgroup\fonteauteurs\bgroup Digne\egroup\egroup{} et Jean
  \bgroup\fonteauteurs\bgroup Michel\egroup\egroup{} :
\newblock Fonctions {$L$} des vari\'et\'es de {D}eligne-{L}usztig et descente
  de {S}hintani.
\newblock {\em M\'em. Soc. Math. France (N.S.)}, (20)\string:\penalty500\relax
  iv+144, 1985.

\bibitem[Dri74]{drinfeld-ell}
V.~G. \bgroup\fonteauteurs\bgroup Drinfel{\cprime}d\egroup\egroup{} :
\newblock Elliptic modules.
\newblock {\em Mat. Sb. (N.S.)}, 94(136)\string:\penalty500\relax 594--627,
  656, 1974.

\bibitem[Dri76]{drinfeld-covering}
V.~G. \bgroup\fonteauteurs\bgroup Drinfel{\cprime}d\egroup\egroup{} :
\newblock Coverings of {$p$}-adic symmetric domains.
\newblock {\em Funkcional. Anal. i Prilo\v zen.},
  10(2)\string:\penalty500\relax 29--40, 1976.

\bibitem[Far08]{FGL}
Laurent \bgroup\fonteauteurs\bgroup Fargues\egroup\egroup{} :
\newblock L'isomorphisme entre les tours de {L}ubin-{T}ate et de {D}rinfeld et
  applications cohomologiques.
\newblock \emph{In} {\em L'isomorphisme entre les tours de {L}ubin-{T}ate et de
  {D}rinfeld}, volume 262 de {\em Progr. Math.}, pages 1--325. Birkh\"auser,
  Basel, 2008.

\bibitem[FvdP04]{put}
Jean \bgroup\fonteauteurs\bgroup Fresnel\egroup\egroup{} et Marius van~der
  \bgroup\fonteauteurs\bgroup Put\egroup\egroup{} :
\newblock {\em Rigid analytic geometry and its applications}, volume 218 de
  {\em Progress in Mathematics}.
\newblock Birkh\"auser Boston Inc., Boston, MA, 2004.

\bibitem[Gen]{genestier}
Alain \bgroup\fonteauteurs\bgroup Genestier\egroup\egroup{} :
\newblock Expos\'es \`a {G}\"ottingen et \`a {H}arvard.

\bibitem[Gen96]{genestier-livre}
Alain \bgroup\fonteauteurs\bgroup Genestier\egroup\egroup{} :
\newblock Espaces sym\'etriques de {D}rinfeld.
\newblock {\em Ast\'erisque}, (234)\string:\penalty500\relax ii+124, 1996.

\bibitem[GI63]{GI}
O.~\bgroup\fonteauteurs\bgroup Goldman\egroup\egroup{} et
  N.~\bgroup\fonteauteurs\bgroup Iwahori\egroup\egroup{} :
\newblock The space of $p$-adic norms.
\newblock {\em Acta Math. 109}, pages 137--177, 1963.

\bibitem[Har77]{hartshorne}
Robin \bgroup\fonteauteurs\bgroup Hartshorne\egroup\egroup{} :
\newblock {\em Algebraic geometry}.
\newblock Springer-Verlag, New York, 1977.
\newblock Graduate Texts in Mathematics, No. 52.

\bibitem[Har97]{Harris-Carayol}
Michael \bgroup\fonteauteurs\bgroup Harris\egroup\egroup{} :
\newblock Supercuspidal representations in the cohomology of {D}rinfel\cprime d
  upper half spaces; elaboration of {C}arayol's program.
\newblock {\em Invent. Math.}, 129(1)\string:\penalty500\relax 75--119, 1997.

\bibitem[Hau05]{Hausberger-these}
Thomas \bgroup\fonteauteurs\bgroup Hausberger\egroup\egroup{} :
\newblock Uniformisation des vari\'et\'es de {L}aumon-{R}apoport-{S}tuhler et
  conjecture de {D}rinfeld-{C}arayol.
\newblock {\em Ann. Inst. Fourier (Grenoble)}, 55(4)\string:\penalty500\relax
  1285--1371, 2005.

\bibitem[HT01]{Harris-Taylor}
Michael \bgroup\fonteauteurs\bgroup Harris\egroup\egroup{} et Richard
  \bgroup\fonteauteurs\bgroup Taylor\egroup\egroup{} :
\newblock {\em The geometry and cohomology of some simple {S}himura varieties},
  volume 151 de {\em Annals of Mathematics Studies}.
\newblock Princeton University Press, Princeton, NJ, 2001.

\bibitem[Ito05]{ito}
Tetsushi \bgroup\fonteauteurs\bgroup Ito\egroup\egroup{} :
\newblock Weight-monodromy conjecture for {$p$}-adically uniformized varieties.
\newblock {\em Invent. Math.}, 159(3)\string:\penalty500\relax 607--656, 2005.

\bibitem[Kur80]{Kurihara}
Akira \bgroup\fonteauteurs\bgroup Kurihara\egroup\egroup{} :
\newblock Construction of {$p$}-adic unit balls and the {H}irzebruch
  proportionality.
\newblock {\em Amer. J. Math.}, 102(3)\string:\penalty500\relax 565--648, 1980.

\bibitem[MS10]{MS}
Ralf \bgroup\fonteauteurs\bgroup Meyer\egroup\egroup{} et Maarten
  \bgroup\fonteauteurs\bgroup Solleveld\egroup\egroup{} :
\newblock Resolutions for representations of reductive {$p$}-adic groups via
  their buildings.
\newblock {\em J. Reine Angew. Math.}, 647\string:\penalty500\relax 115--150,
  2010.

\bibitem[Mus78]{Mustafin}
G.~A. \bgroup\fonteauteurs\bgroup Mustafin\egroup\egroup{} :
\newblock Non-{A}rchimedean uniformization.
\newblock {\em Mat. Sb. (N.S.)}, 105(147)(2)\string:\penalty500\relax 207--237,
  287, 1978.

\bibitem[Rap90]{Rapoport}
M.~\bgroup\fonteauteurs\bgroup Rapoport\egroup\egroup{} :
\newblock On the bad reduction of {S}himura varieties.
\newblock \emph{In} {\em Automorphic forms, {S}himura varieties, and
  {$L$}-functions, {V}ol.\ {II} ({A}nn {A}rbor, {MI}, 1988)}, volume~11 de {\em
  Perspect. Math.}, pages 253--321. Academic Press, Boston, MA, 1990.

\bibitem[Ray74]{Raynaud}
Michel \bgroup\fonteauteurs\bgroup Raynaud\egroup\egroup{} :
\newblock Sch\'emas en groupes de type {$(p,\dots, p)$}.
\newblock {\em Bull. Soc. Math. France}, 102\string:\penalty500\relax 241--280,
  1974.

\bibitem[RZ96]{RZ}
M.~\bgroup\fonteauteurs\bgroup Rapoport\egroup\egroup{} et Th.
  \bgroup\fonteauteurs\bgroup Zink\egroup\egroup{} :
\newblock {\em Period spaces for {$p$}-divisible groups}, volume 141 de {\em
  Annals of Mathematics Studies}.
\newblock Princeton University Press, Princeton, NJ, 1996.

\bibitem[SGA1]{SGA1}
{\em Rev\^etements \'etales et groupe fondamental}.
\newblock Lecture Notes in Mathematics, Vol. 224. Springer-Verlag, Berlin,
  1971.
\newblock S{\'e}minaire de g{\'e}om{\'e}trie alg{\'e}brique du Bois Marie
  1960--61.

\bibitem[SGA7-1]{SGA7-1}
{\em Groupes de monodromie en g\'eom\'etrie alg\'ebrique. {I}}.
\newblock Lecture Notes in Mathematics, Vol. 288. Springer-Verlag, Berlin,
  1972.
\newblock S{\'e}minaire de G{\'e}om{\'e}trie Alg{\'e}brique du Bois-Marie
  1967--1969 .

\bibitem[SGA4-3]{SGA4-3}
{\em Th\'eorie des topos et cohomologie \'etale des sch\'emas. {T}ome 3}.
\newblock Lecture Notes in Mathematics, Vol. 305. Springer-Verlag, Berlin,
  1973.
\newblock S{\'e}minaire de G{\'e}om{\'e}trie Alg{\'e}brique du Bois-Marie
  1963--1964 .

\bibitem[SS93]{SS-crelle}
P.~\bgroup\fonteauteurs\bgroup Schneider\egroup\egroup{} et
  U.~\bgroup\fonteauteurs\bgroup Stuhler\egroup\egroup{} :
\newblock Resolutions for smooth representations of the general linear group
  over a local field.
\newblock {\em J. Reine Angew. Math.}, 436\string:\penalty500\relax 19--32,
  1993.

\bibitem[Str08]{Strauch}
Matthias \bgroup\fonteauteurs\bgroup Strauch\egroup\egroup{} :
\newblock Deformation spaces of one-dimensional formal modules and their
  cohomology.
\newblock {\em Adv. Math.}, 217(3)\string:\penalty500\relax 889--951, 2008.

\bibitem[Tei90]{teitelbaum}
Jeremy \bgroup\fonteauteurs\bgroup Teitelbaum\egroup\egroup{} :
\newblock Geometry of an \'etale covering of the {$p$}-adic upper half plane.
\newblock {\em Ann. Inst. Fourier (Grenoble)}, 40(1)\string:\penalty500\relax
  68--78, 1990.



\bibitem[Wan13]{wang_DL}
Haoran \bgroup\fonteauteurs\bgroup Wang\egroup\egroup{} :
\newblock {S}ur la cohomologie des compactifications de vari\'et\'es de
  {D}eligne-{L}usztig.
\newblock {\em \`A para\^itre dans Ann. Inst. Fourier (Grenoble)},
\newblock {ar{X}iv:1310.7259, 2013.}

\bibitem[Wan14]{Wang-Sigma2}
Haoran \bgroup\fonteauteurs\bgroup Wang\egroup\egroup{} :
\newblock L'espace sym\'etrique de {D}rinfeld et correspondance de {L}anglands locale {II},
\newblock {arXiv:1402.1965}, 2014.

\bibitem[Yos10]{Yoshida}
Teruyoshi \bgroup\fonteauteurs\bgroup Yoshida\egroup\egroup{} :
\newblock On non-abelian {L}ubin-{T}ate theory via vanishing cycles.
\newblock \emph{In} {\em Algebraic and arithmetic structures of moduli spaces
  ({S}apporo 2007)}, volume~58 de {\em Adv. Stud. Pure Math.}, pages 361--402.
  Math. Soc. Japan, Tokyo, 2010.

\bibitem[Zhe08]{zheng}
Weizhe \bgroup\fonteauteurs\bgroup Zheng\egroup\egroup{} :
\newblock Sur la cohomologie des faisceaux {$l$}-adiques entiers sur les corps
  locaux.
\newblock {\em Bull. Soc. Math. France}, 136(3)\string:\penalty500\relax
  465--503, 2008.

\end{thebibliography}
\end{document}